\documentclass[a4paper,12pt]{article}%

\usepackage[latin1]  {inputenc}%
\usepackage[T1]      {fontenc }%
\usepackage          {ae      }%
\usepackage          {latexsym}%
\usepackage          {amsmath }%
\usepackage          {amsfonts}%
\usepackage          {amssymb }%
\usepackage          {amsthm  }%
\usepackage[dvips]   {graphicx}%
\usepackage          {wrapfig }%
\usepackage          {floatflt}%
\usepackage          {fancyhdr}%
\usepackage          {a4      }%
\usepackage          {psfrag  }%
\usepackage          {url     }%
\usepackage          {paralist}%
\usepackage          {color   }%
\usepackage          {calc    }%
\usepackage[small,it]{caption }%

\usepackage[totalwidth=.7\paperwidth,
            totalheight=.85\paperheight,
            headsep=.01\paperheight,
            hmarginratio=2:3,
            vmarginratio=2:3,
            includeheadfoot,
            bindingoffset=5mm,
            dvips]{geometry}%

\DeclareMathOperator{\flag }{flag }%
\DeclareMathOperator{\conv }{conv }%
\DeclareMathOperator{\rs}{\mathcal R}%
\DeclareMathOperator{\interior}{int}%

\newcommand{\R}{\mathbb R}%
\newcommand{\N}{\mathbb N}%
\newcommand{\Z}{\mathbb Z}%

\newcommand{\skp}[2]{\langle #1,#2\rangle}%
\newcommand{\card}[1]{|#1|}%
\newcommand{\boundary}{\partial}

\allowdisplaybreaks[4]%

\theoremstyle{plain}%
\newtheorem{theorem}{Theorem}[section]%
\newtheorem{cor}[theorem]{Corollary}%
\newtheorem{lemma}[theorem]{Lemma}%
\newtheorem{defi}[theorem]{Definition}%

\theoremstyle{definition}%
\newtheorem{remark}[theorem]{Remark}%

\numberwithin{theorem}{section}%
\numberwithin{figure}{section}%
\numberwithin{table}{section}%

\fancypagestyle{plain}{\fancyhead{}\fancyfoot{}\fancyfoot[L]{\footnotesize
    $^*$Research supported by the Deutsche Forschungsgemeinschaft
    within the European graduate program ``Combinatorics, Geometry,
    and Computation'' (No. GRK 588/2)}
                       }%

\makeatletter%
\renewcommand{\subsubsection}{%
  \@startsection{subsubsection}{3}{0mm}{-0.25ex \@plus.3ex \@minus.3ex}{-1em}%
  {\normalfont\normalsize\scshape}%
  }%
\renewcommand{\subsection}{%
  \@startsection{subsection}{2}{0mm}{-1\baselineskip \@plus.2ex \@minus.2ex}{0.3\baselineskip}%
  {\normalfont\normalsize\itshape}%
  }%
\renewcommand{\section}{%
  \@startsection{section}{1}{0mm}{-1\baselineskip \@plus.2ex \@minus.2ex}{0.3\baselineskip}%
  {\normalfont\scshape\centering}%
  }%
\renewcommand\paragraph{%
  \@startsection{paragraph}{4}{0mm}{-0.25ex \@plus.3ex \@minus.3ex}{-1em}%
  {\normalfont\normalsize\bfseries\itshape}%
  }%
\makeatother%

\begin{document}
\thispagestyle{plain}

\vspace*{2\baselineskip}
{
{\large\bfseries\sc New Polytopes from Products}

\vspace*{\baselineskip}
{\sc\small Andreas Paffenholz$^*$}

\medskip
\footnotesize 
Inst.\ Mathematics, MA 6-2, TU Berlin, D-10623 Berlin, Germany

\footnotesize
\url{paffenholz@math.tu-berlin.de}

\medskip
\scriptsize December 2004, revised October 2005

\vspace*{.2\baselineskip}
}

\begin{quote}
  \footnotesize {\bf  Abstract.   }  We construct  a new $2$-parameter
  family    $E_{mn}$,  $m,n\ge    3$,  of  self-dual    $2$-simple and
  $2$-simplicial  $4$-polytopes, with flexible geometric realisations. 
  $E_{44}$ is the $24$-cell.   For  large $m,n$ the $f$-vectors   have
  ``fatness'' close to $6$. 
  
  The $E_t$-construc\-tion   of Paffenholz    and Ziegler   applied to
  products of polygons yields  cellular spheres with the combinatorial
  structure of $E_{mn}$.  Here we prove polytopality of these spheres.
  More   generally,  we construct polytopal   realisations for spheres
  obtained from the  $E_t$-construction     applied to  products    of
  polytopes in any dimension $d\ge3$, if  these polytopes satisfy some
  consistency conditions.
  
  We show  that the  projective realisation  space  of $E_{33}$ is  at
  least   nine   dimensional and  that  of    $E_{44}$ at   least four
  dimensional.  This   proves that the  $24$-cell  is not projectively
  unique.   All   $E_{mn}$  for  relatively   prime  $m,n\ge  5$  have
  automorphisms of  their face    lattice  not induced  by  an  affine
  transformation   of    any  geometric    realisation.   The    group
  $\Z_m\times\Z_n$  generated  by rotations in  the two  polygons is a
  subgroup of   the  automorphisms of  the   face lattice of $E_{mn}$. 
  However,  there are only five pairs  $(m,n)$ for which this subgroup
  is geometrically realisable.
\end{quote}  


\section*{Introduction}

In 2003, Eppstein, Kuperberg, and  Ziegler introduced a new method for
the   construction of  $2$-simple  and   $2$-simplicial  $4$-polytopes
\cite{EKZ02}.  This  was subsequently extended to arbitrary dimensions
and to spheres  and lattices  by Paffenholz and   Ziegler~\cite{PZ03}.
The  construction  produces  PL  $(d-1)$-spheres from $d$-polytopes by
subdividing and combining faces of the polytope  in a certain way.  It
is unknown whether  these spheres are  polytopal in general.  However,
Paffenholz and Ziegler~\cite{PZ03} list several  series of examples in
which they have a polytopal realisation.

Here  we provide sufficient  conditions   for the polytopality of  the
spheres   that we  obtain  when the    con\-struc\-tion is applied  to
products of  two polytopes.   We  present examples  of $d$-dimensional
products for which these conditions  are  satisfied, for all  $d\ge3$.
Our main interest is in the application to products $C_m\times C_n$ of
two polygons with $m$ and $n$ vertices.   We prove that these products
satisfy  our  conditions  for     all  $m,n\ge 3$,  resulting    in  a
two-parameter  family  $E_{mn}$  of   self-dual,  $2$-simplicial   and
$2$-simple polytopes.  All  these polytopes have a large combinatorial
symmetry   group and  only   three  different  combinatorial types  of
vertices and facets.

The   underlying CW spheres in  the  special case $m=n$ were described
earlier  by G{\'e}vay \cite{Gevay04}  and Bokowski  \cite{Bokowski04}.
G{\'e}vay  also   considered symmetry  properties  of  these  spheres.
Polytopality  for  $\tfrac{1}{m}+ \tfrac{1}{n}\ge  \tfrac{1}{2}$  is a
consequence of a theorem of Santos \cite[Rem.~13]{Santos00}.

There are two different notions of symmetry for a polytope:
\begin{inparaenum}[(1)]
\item automorphisms of  the face  lattice (combinatorial  symmetries),
  and
\item transformations  that set-wise preserve  a geometric realisation
  of the polytope (geometric symmetries).
\end{inparaenum}
Any  geometric symmetry    preserves  incidences and  thus  induces  a
combinatorial symmetry. However, the opposite  implication is not true
in general,  i.e.\ not all combinatorial  symmetries of a polytope can
always be realised geometrically in some  realisation of the polytope.
Mani  \cite{Mani71}  and Perles \cite[p.~120]{Gruenbaum67} proved that
all $3$-polytopes, and all  $d$-polytopes with at most $d+3$ vertices,
have  a geometric realisation    whose  geometric symmetry group    is
isomorphic to the combinatorial  one,  while Bo\-kow\-ski, Ewald,  and
Kleinschmidt \cite{BEK84}   presented a $4$-polytope  on $10$ vertices
having a combinatorial symmetry not induced by a geometric one.

Here we prove that all polytopes $E_{mn}$ for relatively prime $m,n\ge
5$   have geometrically     non-realisable  combinatorial  symmetries. 
Furthermore,  the   combinatorial symmetry group   of  $E_{mn}$ always
contains the product $\Z_m\times\Z_n$ of two  cyclic groups induced by
a rotation  of the vertices in  the two  polygons.  However, there are
only five pairs $(m,n)$ in which the  geometric symmetry group of some
realisation has a subgroup inducing these combinatorial symmetries.

The polytope $E_{44}$   is combinatorially equivalent to  the 24-cell,
and applying the  $E$-construction to the  product of two unit squares
produces   its   regular   realisation.    However, our   polytopality
conditions for the $E$-construction of  products  allow for much  more
flexibility.   For the smallest  instance  $E_{33}$  we work out   all
degrees  of  freedom that our conditions   permit and give an explicit
construction of all possible  such realisations.  This will prove that
the projective  realisation space $\rs_{proj}(E_{33})$ of  $E_{33}$ is
at least  nine dimensional.   For  the $24$-cell  we present a  simple
$4$-parameter family of realisations showing that $\rs_{proj}(E_{44})$
is   at  least  four dimensional.    In  particular, the  $24$-cell is
projectively not unique (cf. McMullen~\cite{McMullen76}).

Eppstein,    Kuperberg,    and Ziegler~\cite{EKZ02}    introduced  the
``fatness''  $F(P)$  of a $4$-polytope   $P$,   which is roughly   the
quotient of the number  of edges and ridges by  the number of vertices
and  facets.   They construct  an  example  with fatness approximately
$5.073$.  For  large $m,n$  our  polytopes $E_{mn}$ will  have fatness
arbitrarily close to $6$.   However, Ziegler \cite{Ziegler04} recently
constructed a new family of $4$-polytopes from projections of products
of polygons whose fatness approaches $9$.

I  am    grateful  to   J.~Bokowski,   G.~G{\'e}vay,   F.~Santos,  and
G.M.~Ziegler for hints and discussions.  I am grateful to the referees
for suggesting a  simpler statement of Theorem \ref{thm:construction},
its  proof   and the  proof of   Theorem~\ref{thm:relspace},   and for
pointing out  a    gap    in the    construction     of $D(n,r)$    in
Section~\ref{sec:polytopality}.


\section{Polytopes, products, and the $E$-construction}

This section gives a short  introduction to polytopes, their products,
and the  $E$-construction.   See \cite{Ziegler95} and  \cite{PZ03} for
more background.
 
\paragraph{Polytopes.} 

A {\em polytope} $P$  is the convex hull of  a finite set of points in
$\R^n$.  Its {\em  dimension}  $d$ is the   smallest dimension  of  an
affine  subspace containing $P$.   $V(P)$ denotes the  {\em set of all
  vertices} of  a $d$-polytope $P$.   Faces of codimension $1$ and $2$
are  called {\em     facets}  and  {\em  ridges}.      Let $f_S$   for
$S\subset\{0,\ldots,d-1\}$ be the number of increasing chains with one
face of dimension $i$ for each  $i\in S$.  The vector collecting these
numbers  is   called  the  {\em flag   vector}  $\flag(P)$.  The  {\em
  $f$-vector}\  is  the  subset of  $\flag(P)$   corresponding to  the
entries with $\card{S}=1$.  We set $f_d(P):=1$.

\paragraph{Products.} 

For $i=0,1$  let      $P_i$ be $d_i$-polytopes  with   flag    vectors
$\flag(P_i)=(f_S(P_i))_{S\subseteq\{0,\ldots,d_i-1\}}$.    The    {\em
  product} $P_0\times P_1$ is the convex hull of
\begin{align*}
  V(P_0\times P_1)&:=\{ (v,w)  \in \R^{d_0+d_1} \mid v\in V(P_0), w\in
  V(P_1)\}.
\end{align*}
Equivalently,  $P_0\times P_1:=\{(v,w)\in\R^{d_0+d_1}\mid  v\in  P_0,\;
w\in   P_1\}$.  It has dimension  $(d_0+d_1)$   and flag vector $\flag
(P_0\times  P_1):=(f_S(P_0\times     P_1))_{S  \subseteq  \{0, \ldots,
  d_0+d_1-1\}}$ with
\begin{align*}
  f_S(P_0\times P_1)&:=f_{(s_1,s_2,\ldots,s_k)}(P_0\times P_1)\\[.1cm]
  &\phantom{:}=\sum_{u_1+v_1=s_1}\sum_{u_2+v_2=s_2}\ldots\sum_{u_k+v_k=s_k}
  f_{(u_1,u_2,\ldots,u_k)}(P_0)f_{(v_1,v_2,\ldots,v_k)}(P_1).
\end{align*}
In  this   formula  we   set $f_{(t_1,t_2,\ldots,t_k)}\equiv0$  unless
$t_1\le  t_2\le\ldots\le   t_k$ and    define $f_{(t_1,   t_2, \ldots,
  t_{i-1},  t_i,   t_{i+1},  \ldots,  t_k)}  :=f_{(t_1,   t_2, \ldots,
  t_{i-1}, t_{i+1}, \ldots, t_k)}$ if $t_i=t_{i+1}$.

We  have defined here the  {\em geometric (orthogonal)} product as the
convex  hull  of all pairs  of  geometrically given vertices.   A more
general definition   would just   require a   polytope combinatorially
equivalent to this.

\paragraph{$\mathbf E$-construction.} 

For our purposes the $E$-construction of a $d$-polytope $P$, $d\ge 2$,
is  best viewed as  a construction that takes  polytopes  as input and
produces  regular  CW spheres (their ``$E$-sphere'')   from them.  The
original definition in \cite{PZ03} depends on  a parameter $t$ between
$0$  and $d-1$ (the dimension of  ``distinguished'' elements). We omit
this parameter in the notation, as we use only the case $t=d-2$.

\begin{figure}[b]
\begin{minipage}[t]{.3\textwidth}
\centering
\includegraphics[width=.95\textwidth]{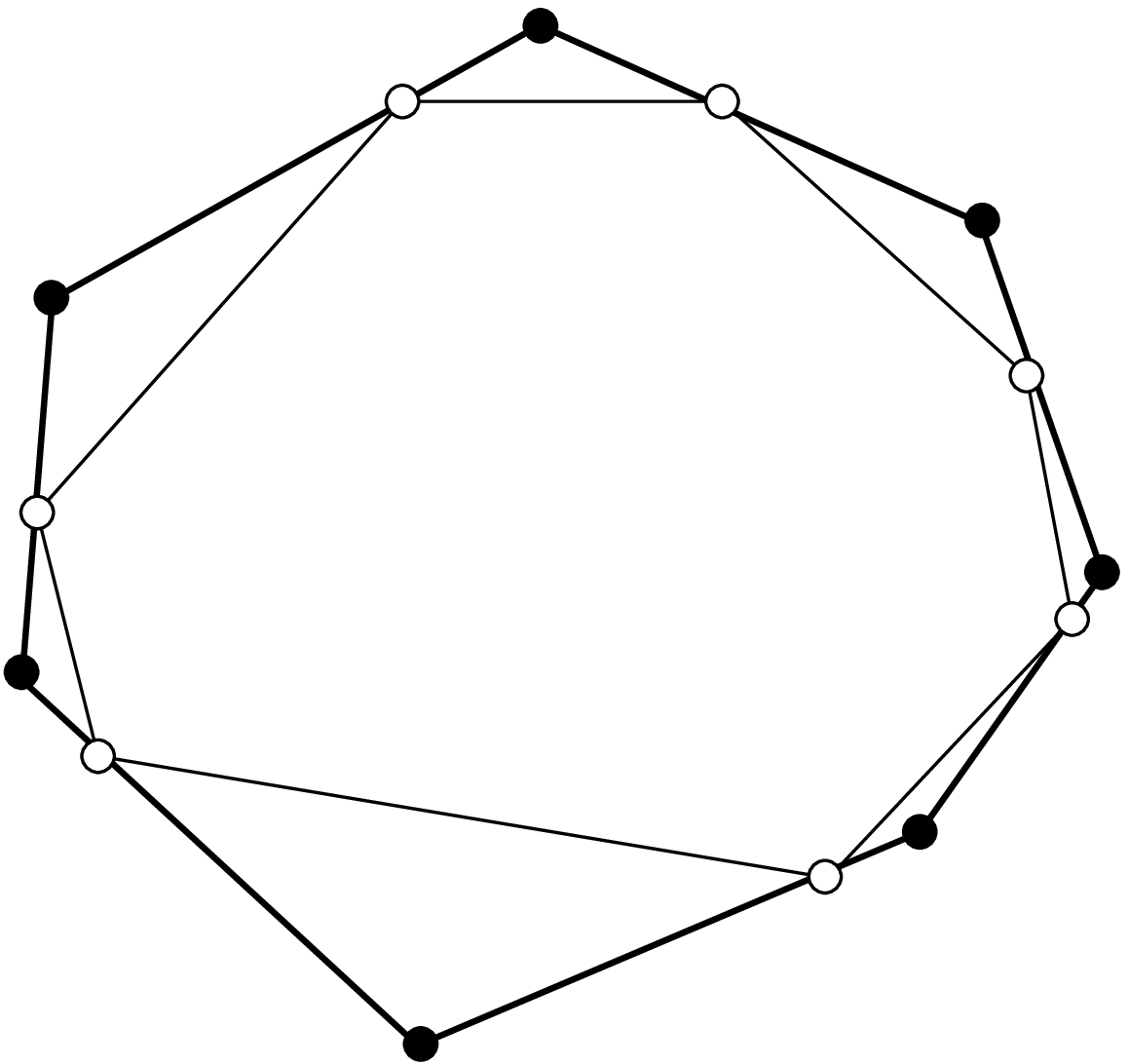}
\caption{The $E$-con\-struc\-tion (thick edges) applied to a polygon 
  (thin edges).}
\label{fig:Etpolygon}
\end{minipage}
\hspace{.02\textwidth}
\begin{minipage}[t]{.6\textwidth}
\centering
\includegraphics[width=.45\textwidth,bb=117 56 508 394,clip]{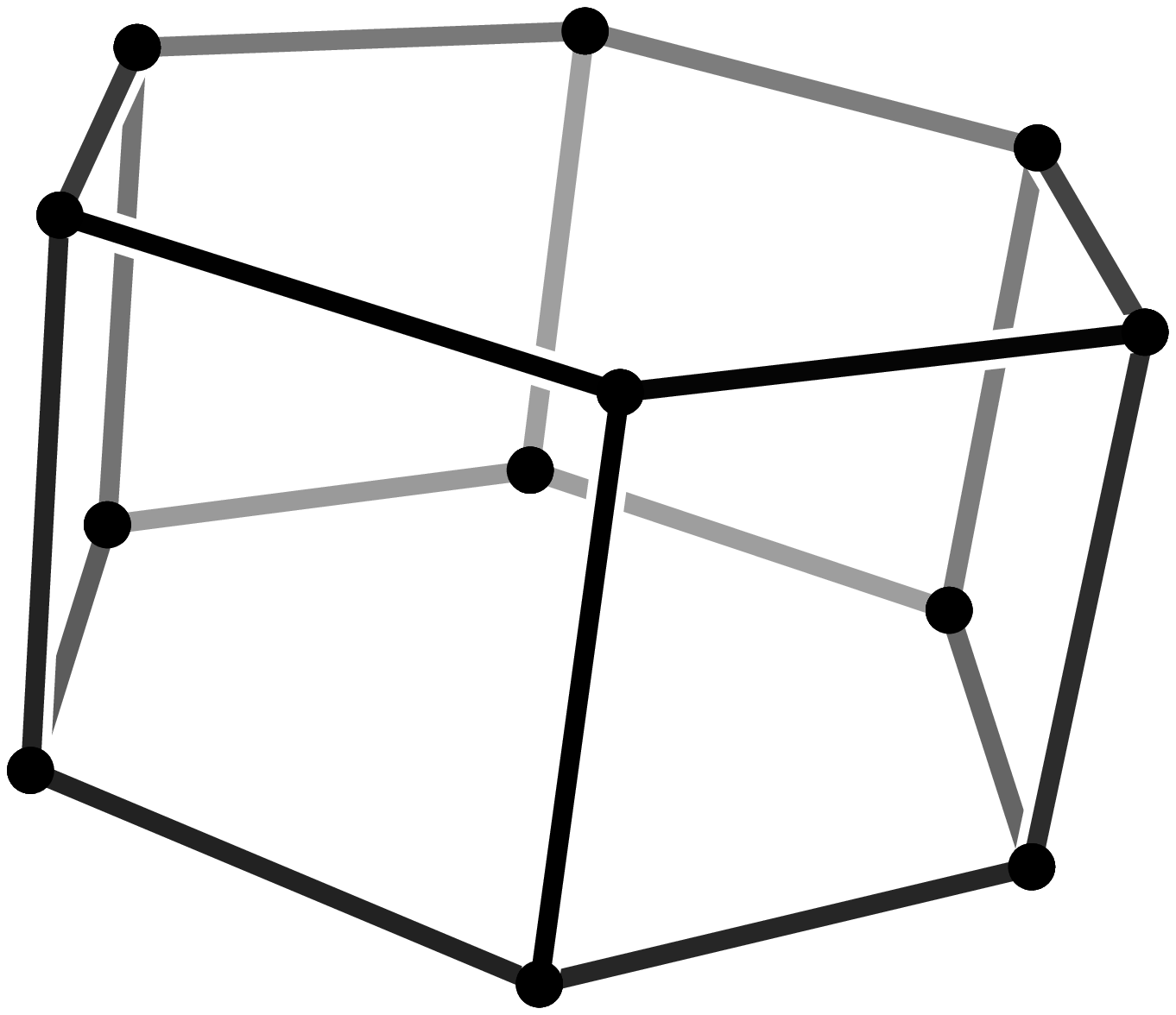}
\hspace{.05\textwidth}
\includegraphics[width=.45\textwidth,bb=30 30 574 570,clip]{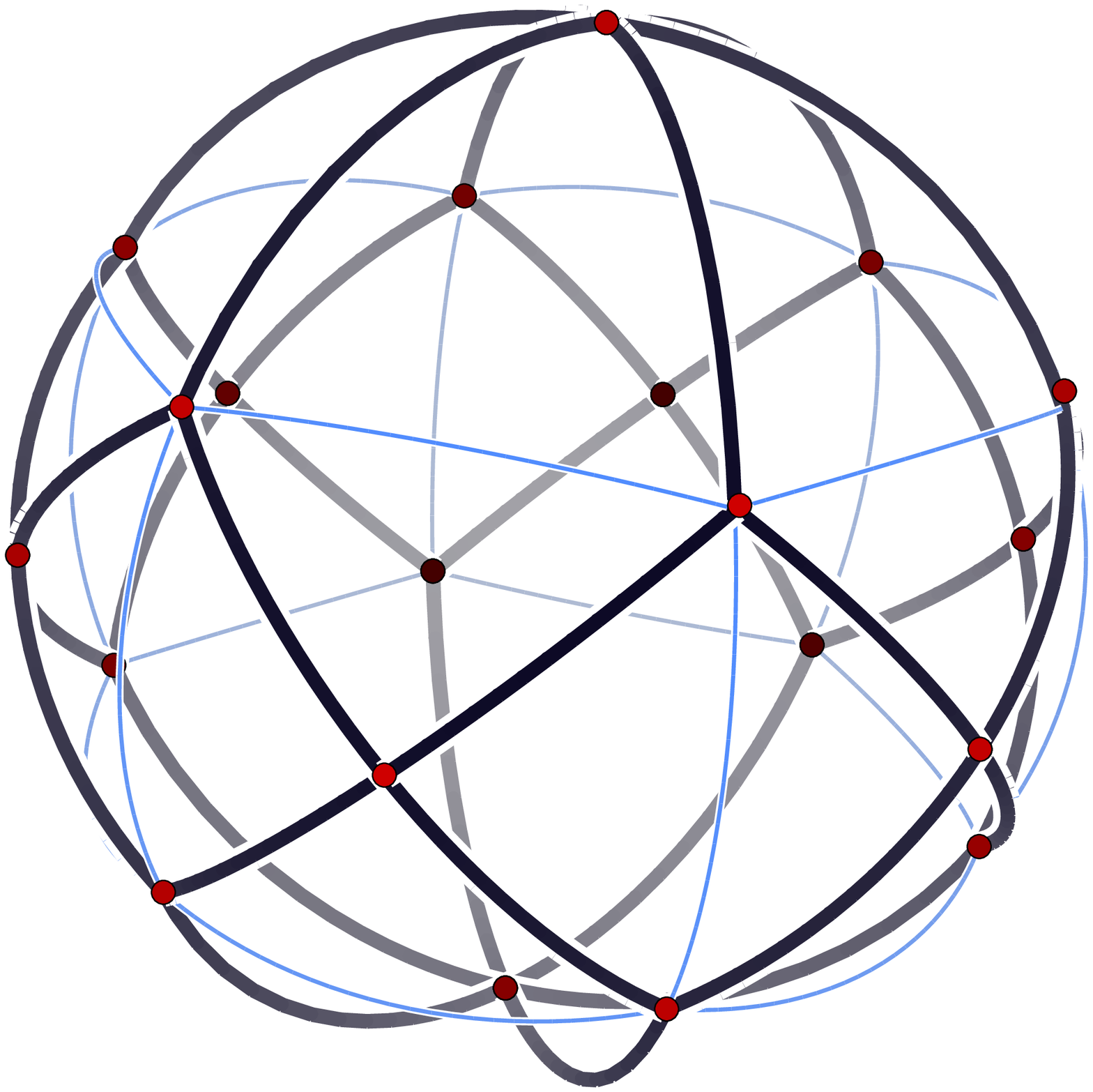}
\caption{A polytope and its $E$-sphere (in bold, the polytope is drawn
  thin to show the old ridges).}
\label{fig:Et3D}
\end{minipage}
\end{figure}

Here is    the  construction.    Let $P$ be    a   $d$-polytope.   The
$E$-construction assigns to $P$ a  $CW$-sphere $E(P)$ by the following
two steps:
\begin{compactenum}[(1)]
\item {\it Stellarly subdivide all facets of the polytope $P$,}
\item {\it  and  merge facets of the  subdivision  sharing a  ridge of
    $P$.}
\end{compactenum}
Each facet of the subdivision contains precisely one such ridge, so we
merge pairs of facets of  the subdivision.  Thus, combinatorially  the
facets  of  $E(P)$   are bipyramids  over    the ridges  of~$P$.   See
Figure~\ref{fig:Etpolygon}    and    Figure~\ref{fig:Et3D}   for     a
two-dimensional and a  three-dimensional example of this construction.
In  dimensions $d\ge  3$ all   vertices of the  original polytope  are
preserved, while for $d=2$  they lie on  the new edges.  A more formal
definition   given   on  the    level    of  face   lattices  is    in
\cite[Def.~1.2]{PZ03}.  For any  polytope  $P$ the sphere $E(P)$  is a
piecewise linear CW  sphere \cite[Thm.  2.1]{PZ03}.  If  these spheres
are polytopal  then   we   call  the   resulting  polytope   the  {\em
  $E$-polytope}     of  $P$.   This  is    e.g.\   the  case  for  all
dual-to-stacked $4$-polytopes \cite[Sect.~3]{PZ03}.  The $f$-vector of
$E(P)$ is given by
\begin{align*}
  f_k(E(P)):=
  \begin{cases}
    f_{d-2}(P)\qquad& k=d-1\\
    f_{d-3,d-1}(P)& k=d-2\\
    f_k(P)+f_{k-1,d-1}(P)\quad& \text{otherwise,}
  \end{cases}
\end{align*}
where we set $f_{-1j}:=f_j$.

\begin{figure}[b]
\centering
\includegraphics[width=.375\textwidth]{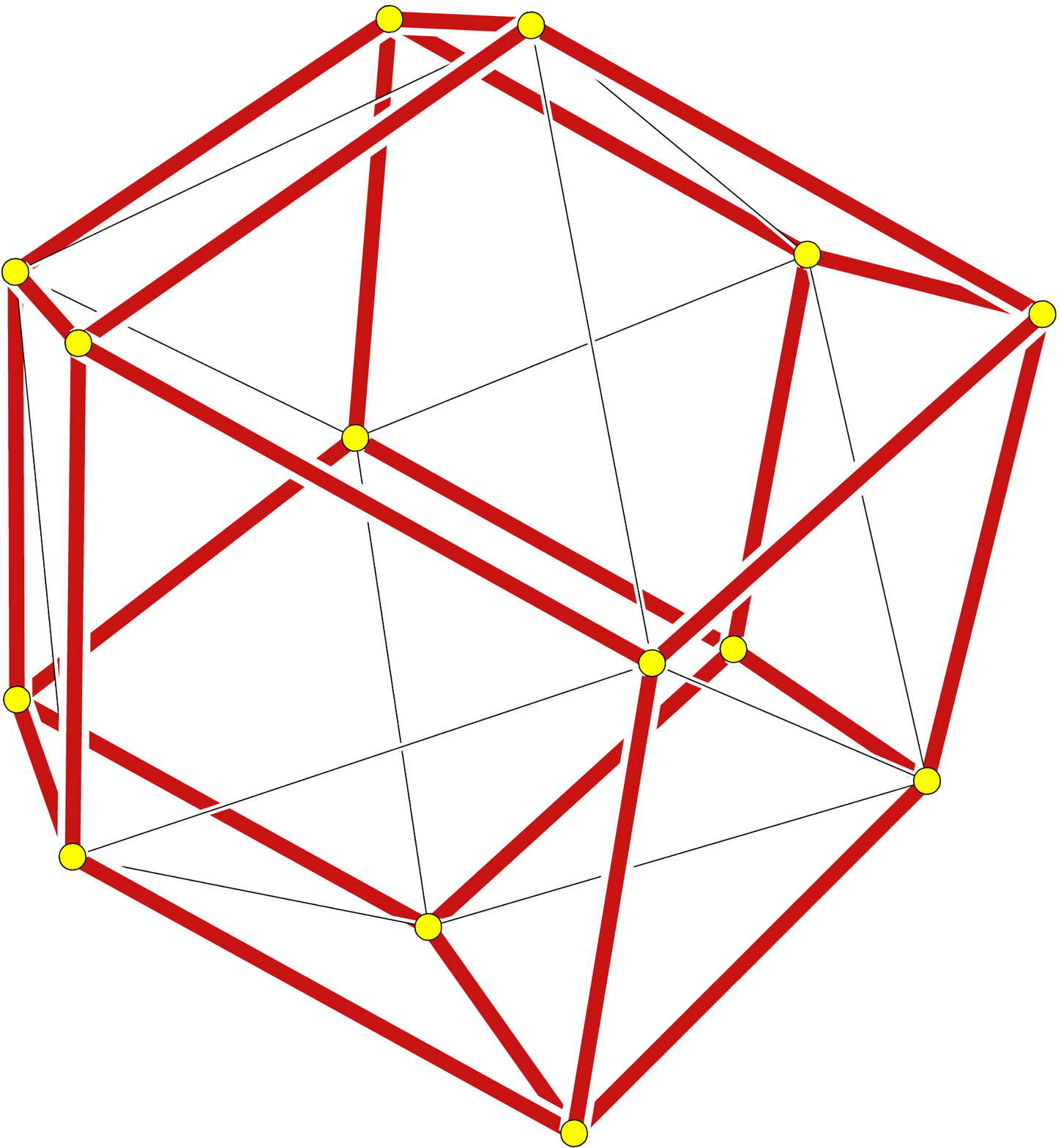}
\hspace{.1\textwidth}
\includegraphics[width=.375\textwidth]{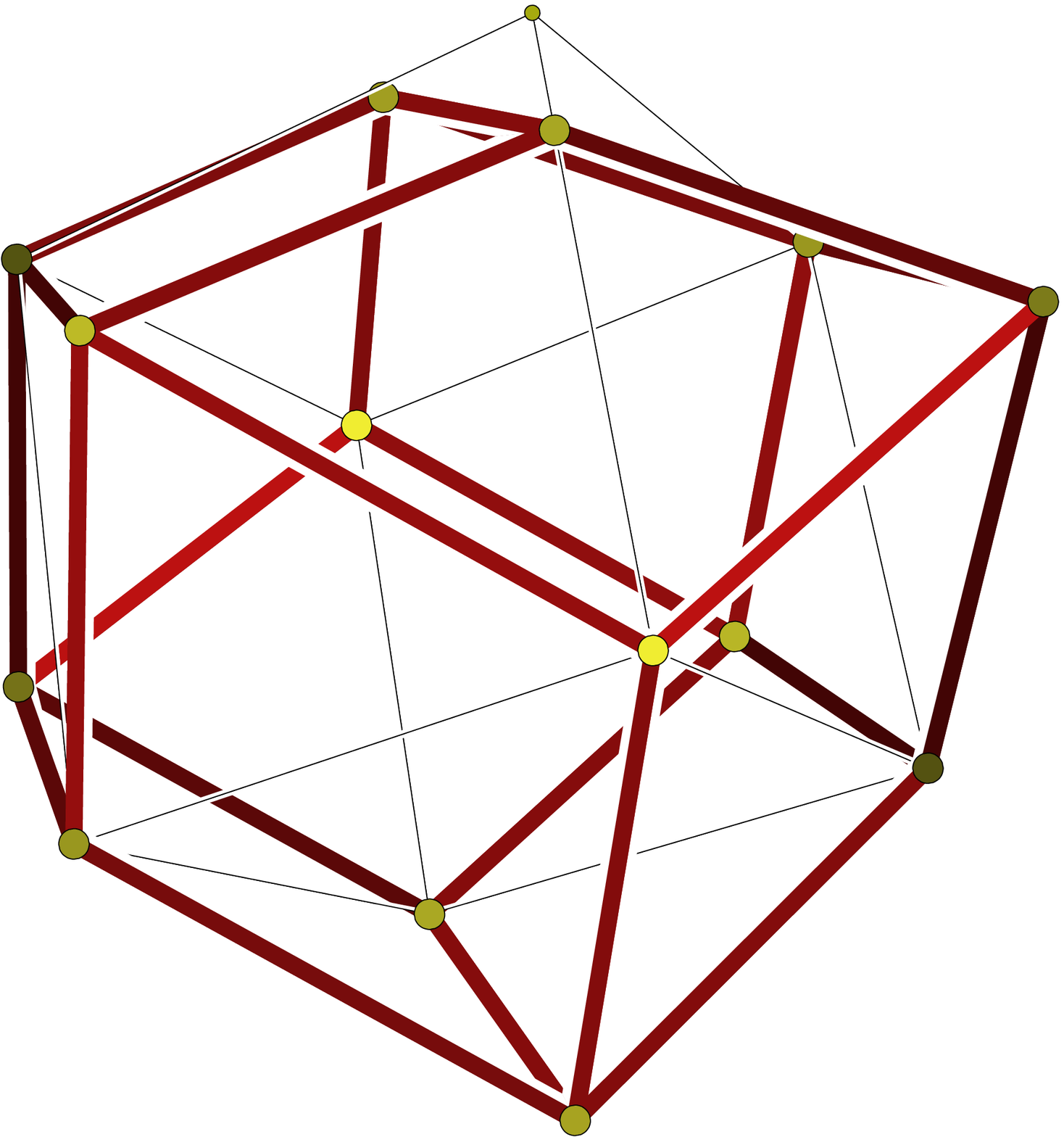}
\caption{The left realisation of $E(C)$ of the unit cube 
  $C$ is  vertex-preserving, the right is not:  observe the top vertex
  of the cube (and there is no cube for which it is).}
\label{fig:3Dcube}
\end{figure}
In the above  definition the $E$-polytope  of  some polytope $P$  just
denotes some polytope  being combinatorially equivalent to  the sphere
obtained from $P$ via the $E$-construction.   In the following we need
a stricter version of the connection between $P$ and its $E$-polytope.
\begin{defi}[vertex-preserving]
  A polytopal realisation of $E(P)$ for a given geometrically realised
  polytope $P$ is called   {\em vertex-preserving} if it  is  obtained
  from the  realisation  of $P$  by  placing new  vertices beyond  the
  facets    of  $P$   and  taking  the     convex   hull. See   Figure
  \ref{fig:3Dcube} for an example.
\end{defi}

\begin{remark}
  For illustrations we will sometimes  also apply the  $E$-construction
  to a $1$-polytope  $S$ (i.e.\ a segment).  In   this case $E(S)$  is
  defined to be a segment containing $S$ in its interior. 
\end{remark}

\paragraph{Polygons.}

We denote a (convex) polygon with $m$  vertices $v_0, \ldots, v_{m-1}$
by  $C_m$.     We  usually  assume that   the   vertices are  numbered
consecutively and take indices modulo $m$. 

By $E_{mn}$ we  denote the result  of the $E$-construction  applied to
the product $C_m\times  C_n$ of an  $m$-gon and an  $n$-gon. This is a
$4$-dimensional $2$-simplicial and  $2$-simple $CW$ sphere.   The flag
vectors of $C_m\times C_n$ and $E_{mn}$ are:
\begin{align}
\flag(C_m\!\times\! C_n)&\!=\!(mn,2mn,mn+m+n,m+n;4mn)\notag\\
\flag(E_{mn})&\!=\!
(mn\!+\!m\!\!+\!n, 6mn,6mn,mn\!+\!m\!+\!n;8mn\!+\!2(m\!+\!n)),\label{eq:fvector}
\end{align}
where we have only recorded  the values $(f_0,\ldots,f_3;f_{03})$. All
other  entries  of   the flag vector    follow   from the  generalised
Dehn-Sommerville equations \cite{BB85}.


\section{The $E$-construction of products}\label{sec:Eprod}

Let $P_0$, $P_1$ be two  polytopes of dimensions  $d_0$ and $d_1$.  We
give sufficient    conditions   for the    existence  of  a  polytopal
realisation   of the   sphere $E(P_0\times  P_1)$   obtained from  the
polytope  $P_0\times  P_1$.    If  we  restrict  to  vertex-preserving
realisations,  then   these   conditions  are  also    necessary.  The
conditions are the following:

\medskip

\noindent
\fbox{\label{prod-suff-cond}
  \parbox{\textwidth-.5cm}{
    \begin{compactenum}[(A)]
    \item\label{conditionA} There exist vertex-preserving realisations
      of $E(P_0)$ and $E(P_1)$.
    \item\label{conditionB}  For  $i=0,1$ let $T_i:=V(E(P_i))\setminus
      V(P_i)$.  There are two maps
      \begin{align*}
        \beta_i:T_i\rightarrow \interior(P_{1-i})
      \end{align*}
      such that for  any $(v_0,v_1)\in T_0\times  T_1$ the fraction of
      the   segment  $|v_0,\beta_1(v_1)|$  outside  $P_0$  equals  the
      fraction of the segment $|v_1,\beta_0(v_0)|$ inside $P_1$.
    \end{compactenum}
  }
}

\begin{figure}[b] 
  \psfrag{P1}[l][l]{{$\scriptstyle P_1$}}
  \psfrag{P0}{\colorbox[rgb]{1,1,1}{$\scriptstyle P_0$}}
  \psfrag{EP1}{{$\scriptstyle E(P_1)$}}
  \psfrag{EP0}{\colorbox[rgb]{1,1,1}{$\scriptstyle E(P_0)$}}
  \psfrag{v01}[r][r]{\colorbox[rgb]{1,1,1}{$\scriptstyle v_0^1$}}
  \psfrag{v02}[l][l]{\colorbox[rgb]{1,1,1}{$\scriptstyle v_0^2$}}
  \psfrag{v03}[rb][rb]{\colorbox[rgb]{1,1,1}{$\scriptstyle v_0^3$}}
  \psfrag{v11}[r][r]{\colorbox[rgb]{1,1,1}{$\scriptstyle v_1^1$}}
  \psfrag{v12}[l][l]{\colorbox[rgb]{1,1,1}{$\scriptstyle v_1^2$}}
  \psfrag{v13}[l][l]{\colorbox[rgb]{1,1,1}{$\scriptstyle v_1^3$}}
  \psfrag{T1v1}[bl][bl]{\colorbox[rgb]{1,1,1}{$\scriptstyle \beta_1(v_1^1)$}}
  \psfrag{T1v2}[bl][bl]{\colorbox[rgb]{1,1,1}{$\scriptstyle \beta_1(v_1^2)$}}
  \psfrag{T1v0}[b][b]{\colorbox[rgb]{1,1,1}{$\scriptstyle \beta_1(v_1^3)$}}
  \psfrag{T0v1}{\colorbox[rgb]{1,1,1}{$\scriptstyle \beta_0(v_0^1)$}}
  \psfrag{T0v2}[rb][rb]{\colorbox[rgb]{1,1,1}{$\scriptstyle \beta_0(v_0^2)$}}
  \psfrag{T0v3}[b][b]{\colorbox[rgb]{1,1,1}{$\scriptstyle \beta_0(v_0^3)$}}
  \psfrag{q0}{\colorbox[rgb]{1,1,1}{$\scriptstyle q_0$}}
  \psfrag{q1}{\colorbox[rgb]{1,1,1}{$\scriptstyle q_1$}}
  \begin{minipage}{7cm}
    \includegraphics[height=.3\textheight]{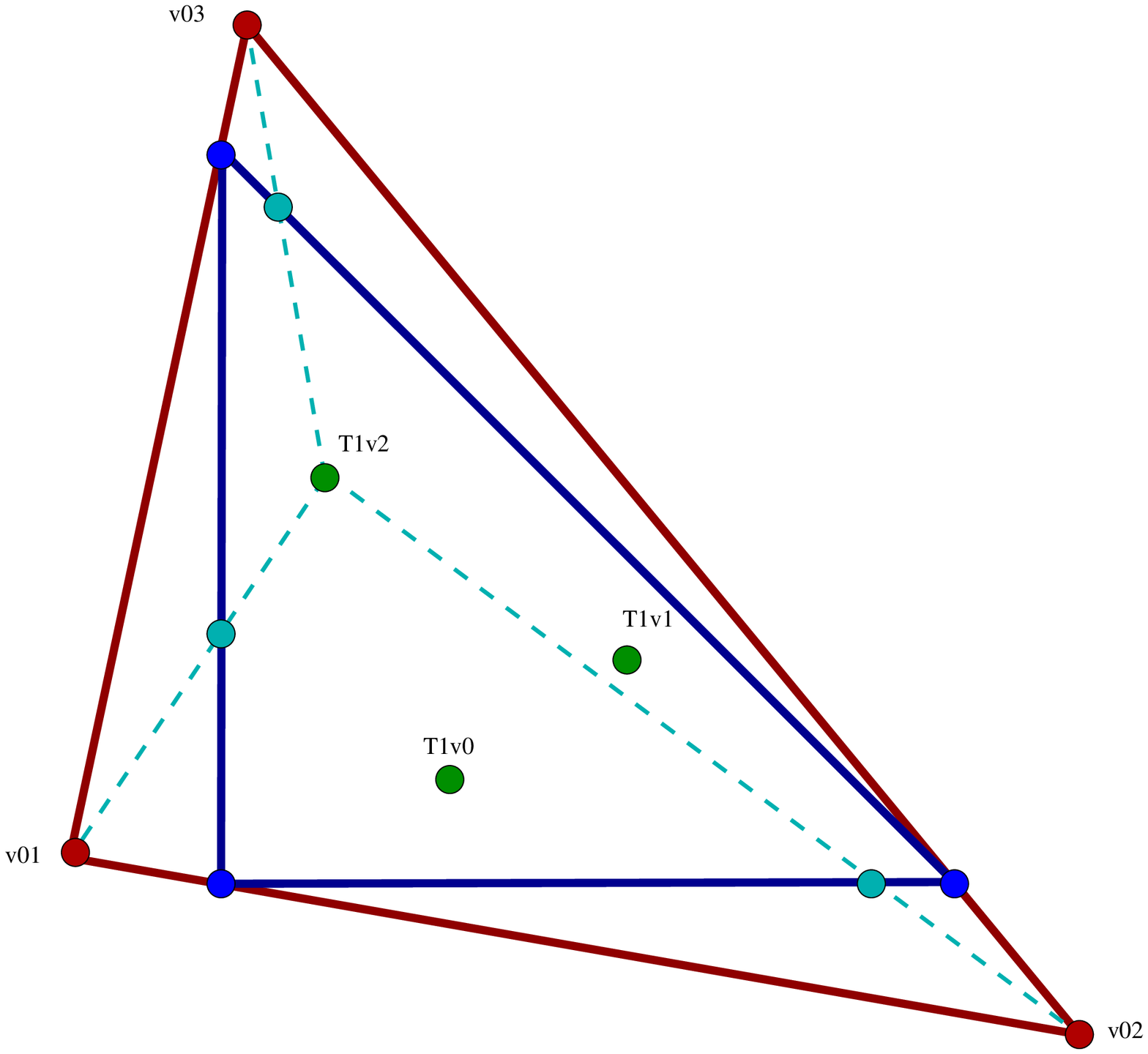}
  \end{minipage} 
  \begin{minipage}{7cm}
    \includegraphics[height=.3\textheight]{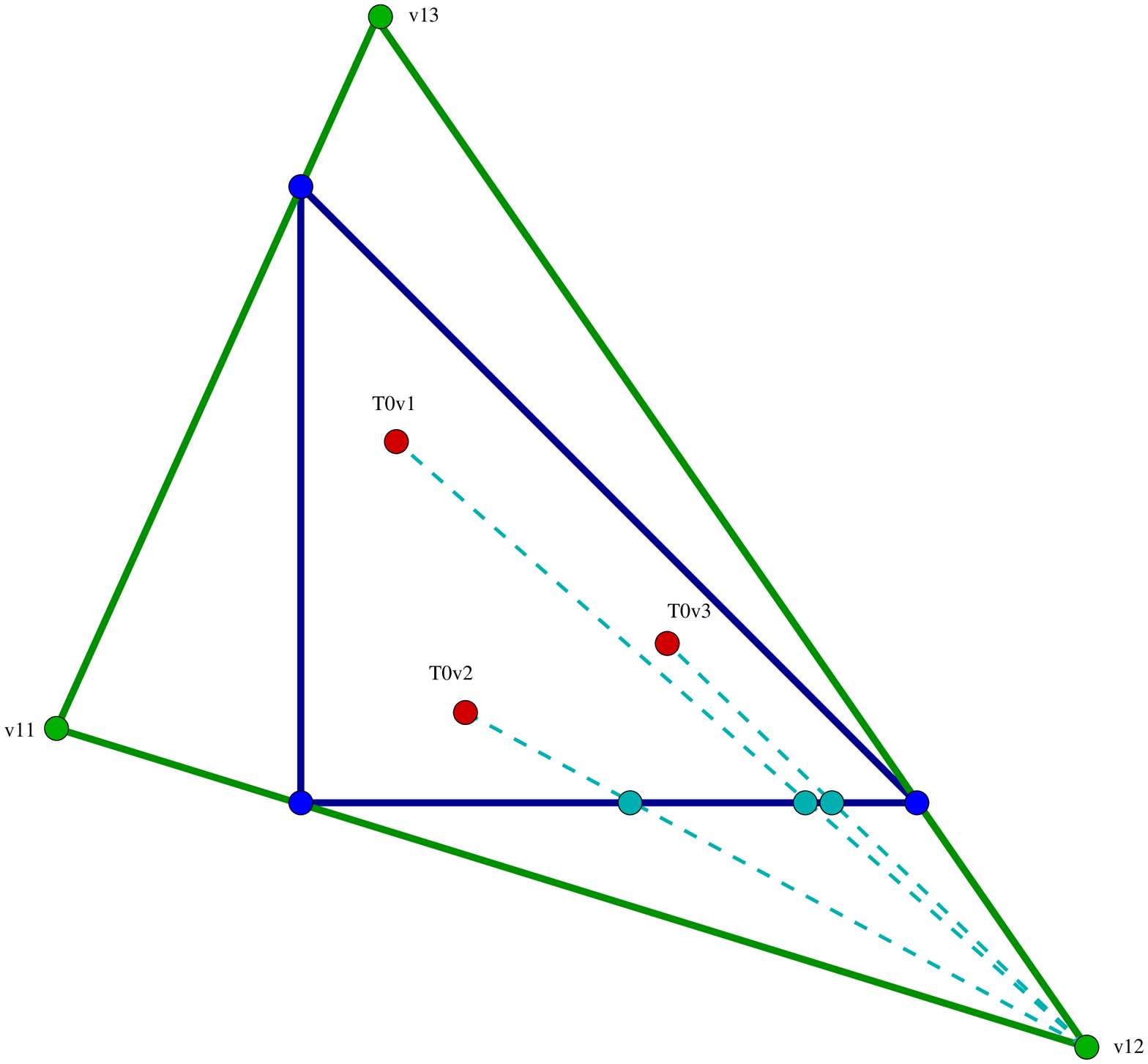}
  \end{minipage}
  \caption{Realising the product of two triangles}
  \label{fig:ProdTriangles}
\end{figure}

\medskip
\begin{theorem}\label{thm:construction}
  Let $P_0$, $P_1$ be a pair of polytopes with $\dim(P_0\times P_1)\ge
  3$        that          satisfies   {\rm(\ref{conditionA})}      and
  {\rm(\ref{conditionB})}.   Let $S\subset\R^{d_0}\times\R^{d_1}$   be
  the point set containing the following points:
  \begin{compactenum}[\rm (a)]
  \item all pairs $(p_0,p_1)$ for $p_0\in V(P_0)$, $p_1\in V(P_1)$,
  \item all pairs $(v_0,\beta_0(v_0))$ for $v_0\in T_0$,
  \item all pairs $(\beta_1(v_1),v_1)$ for $v_1\in T_1$.
  \end{compactenum}
  Then  $\conv(S)$  is a   vertex-preserving polytopal  realisation of
  $E(P_0\times  P_1)$.     Moreover,     for    the    existence    of
  vertex-preserving   realisations  of   $E(P_0\times P_1)$  the   two
  conditions {\rm (\ref{conditionA})} and {\rm (\ref{conditionB})} are
  both necessary and sufficient.
\end{theorem}
\noindent
See  Figure \ref{fig:ProdTriangles}  for an  example  of two triangles
satisfying    {\rm(\ref{conditionA})}     and {\rm(\ref{conditionB})}.
\fboxsep1pt

\begin{proof}
  The proof has two  parts.  First we  prove the necessity of  the two
  conditions  (\ref{conditionA})  and  (\ref{conditionB}) for   vertex
  preserving realisations and then their sufficiency.

  Let $P_0$  and  $P_1$ be  two   geometrically realised  polytopes of
  dimension   $d_0$   and   $d_1$   with  $d_0+d_1~\ge~3$.     Suppose
  $E(P_0\times  P_1)$      exists   and  is      a   vertex-preserving
  rea\-li\-sa\-tion  of $P_0\times P_1$.   We can split the vertex set
  of $E(P_0\times P_1)$  into the vertex  set of $P_0\times P_1$ and a
  set consisting of one vertex beyond each facet of $P_0\times P_1$.

  Define standard  projections   $\pi_j:   \R^{d_0+d_1}\longrightarrow
  \R^{d_j}$ for $j=0,1$.  By  assumption, the vertex  set of $P_j$  is
  contained in $\pi_j(V(E(P_0\times P_1)))$  for $j=0,1$. We determine
  the images of the other vertices of $E(P_0\times P_1)$ under $\pi_0$
  and $\pi_1$.

  \noindent
  The facets of the product $P_0\times P_1$ are of the form
  \begin{inparaenum}[(1)]
    \item  ``Facet    of $P_0$'' $\times   P_1$ or 
    \item    $P_0\times$  ``Facet of $P_1$''.
  \end{inparaenum}
  Thus, we have two different types of ridges:
  \begin{inparaenum}[(I)]
  \item Those between two adjacent facets of the first or second type,
    and
  \item those between a facet of the first and one of the second type.
  \end{inparaenum}
  We deal with these two cases separately:
  \begin{compactenum}[(I)]
  \item Let $F$ and $F'$ be two adjacent facets  of the first type and
    $v$, $v'$  the two vertices  of $E(P_0\times  P_1)$ beyond $F$ and
    $F'$.  Let $R$ be the ridge between $F$ and $F'$.  The projections
    $\pi_0(F)$ and  $\pi_0(F')$  are  adjacent facets   of  $P_0$ with
    common ridge  $\pi_0(R)$.  $\pi_0(v)$  and $\pi_0(v')$  are points
    beyond  these  facets.  $v$, $v'$ and  $R$  lie on a common (facet
    defining) hyperplane  $H$ of $E(P_0\times P_1)$ in $\R^{d_0+d_1}$.
    So the points $\pi_0(v)$, $\pi_0(v')$ and the ridge $\pi_0(R)$ all
    lie on  the  hyperplane  $\pi_0(H)$ in $\R^{d_0}$   and $\pi_0(H)$
    defines a face of $\pi_0(E(P_0\times P_1))$, which must in fact be
    a facet.  Thus, the convex hull of the  projection of all vertices
    of $E(P_0\times P_1)$   is $E(P_0)$.  Similarly, projecting   with
    $\pi_1$ gives $E(P_1)$.
  \item Let $w_0$ and $w_1$ be two vertices of $E(P_0\times P_1)$, the
    first beyond a facet $G_0\times P_1$, the second beyond $P_0\times
    G_1$,  where $G_0$ and  $G_1$ are facets  of $P_0$ and $P_1$.  Let
    $R=G_0\times G_1$  be the  ridge between  these  two facets.   The
    segment $s$ between $w_0$ and $w_1$ intersects $R$ in a point $q$.
    $\pi_0(q)$ is contained in $G_0$  and  $\pi_1(q)$ is contained  in
    $G_1$.  So $\pi_0(w_1)$ is contained in  the interior of $P_0$ and
    $\pi_1(w_0)$ in   the  interior of  $P_1$.   Projections  preserve
    ratios, so
    \begin{align*}        
      r:=\frac{|w_0q|}{|w_0w_1|}
      =\frac{|\pi_0(w_0)\pi_0(q)|}{|\pi_0(w_0)\pi_0(w_1)|}
      =\frac{|\pi_1(w_0)\pi_1(q)|}{|\pi_1(w_0)\pi_1(w_1)|}.
    \end{align*} 
  \end{compactenum}
  Hence, a vertex-preserving realisation of $E(P_0\times P_1)$ implies
  the    conditions  (\ref{conditionA}) and (\ref{conditionB}).   This
  proves the necessity part of the theorem.

  Now    we    prove    sufficiency   of     (\ref{conditionA})    and
  (\ref{conditionB}).     Suppose      we    have,  according       to
  (\ref{conditionA}) and (\ref{conditionB}), constructed $E(P_0)$  and
  $E(P_1)$ and the maps $\beta_i:T_i\rightarrow\interior(P_{1-i})$ for
  $T_i:=V(E(P_i))\setminus V(P_i)$, $i=0,1$, and  have formed the  set
  $S$ defined in the theorem.  We have to show that  all facets of the
  convex hull  of  $S$ defined thereby  are  bipyramids over ridges of
  $P_0\times P_1$ and that there is precisely one vertex of $S$ beyond
  each facet of  $P_0\times P_1$.  There  are  two different  cases to
  consider:
  \begin{compactenum}[(I)]
  \item Let $R$ be a ridge of $P_0$. Then  $R\times P_1$ is a ridge of
    $P_0\times  P_1$.   Let $F$ and $F'$  be  the  two facets of $P_0$
    adjacent  to $R$ and $p$, $p'$  the vertices of $E(P_0)$ above $F$
    and $F'$ (see Figure~\ref{fig:RidgeXPoly}).   Let $v$ be the facet
    normal of the facet $F_E$ of $E(P_0)$ formed by  $R$, $p$ and $p'$
    and let $l:=\skp{v}{p}$.
    \begin{figure}[b]
      \centering
      \psfrag{R}{$R$} 
      \psfrag{F0}{$F$}
      \psfrag{F1}{$F'$}
      \psfrag{p0}{$p$} 
      \psfrag{p1}{$p'$}
      \psfrag{v}{$v$}
      \includegraphics[width=.25\textwidth]{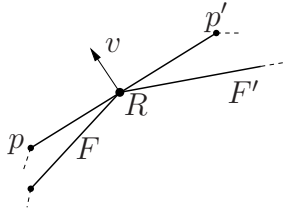}
      \caption{Sufficiency: The case of ``ridge $\times$ polytope''}
      \label{fig:RidgeXPoly}
    \end{figure} 
 
    By construction,  the points $(p,\beta_0(p))$,  $(p',\beta_0(p'))$
    and $(r,q)$,   $r\in V(R)$, $q\in   V(P_1)$ are contained   in the
    hyperplane   $H:=\{x\mid \skp{(v,     \mathbf 0)}{x}=l\}$,   where
    $\mathbf 0$ denotes the $d_1$-dimensional zero vector.  All points
    in the set  $V(E(P_0))\setminus(V(R)\cup\{p,p'\})$ are on the same
    side of the hyperplane defined by the facet  $F_E$.  So all points
    in
    \begin{align*}
      V(E(P_0       \times   P_1))  \setminus(V(R\times  P_1)\cup\{(p,
      \beta_0(p)), (p', \beta_0(p'))\})
    \end{align*}
    are on the same side of the hyperplane $H$ and
    \begin{align*}
      \conv(V(R\times P_1)\cup\{(p,\beta_0(p)), (p',\beta_0(p'))\})
    \end{align*}
    is  a facet of $E(P_0\times P_1)$.    The same argument applies to
    ridges of type $P_0\times R$ for ridges $R$ of $P_1$.
  \item Now consider a ridge of type $F_0\times F_1$ for a facet $F_0$
    of $P_0$ and a  facet $F_1$ of $P_1$.   Let $p_0$ be the vertex of
    $E(P_0)$  beyond  $F_0$ and $p_1$  the   vertex of $E(P_1)$ beyond
    $F_1$.  Let $i_0$ be the intersection point of the segment between
    $p_0$ and  $\beta_1(p_1)$  and  the  facet $F_0$, and   $i_1$  the
    intersection point of the  segment between $p_1$  and $\beta(p_0)$
    and the facet $F_1$. By construction we have
    \begin{align*}
      \frac{|p_0,i_0|}{|p_0,\beta_1(p_1)|}=
      \frac{|\beta_0(p_0),i_1|}{|\beta_0(p_0),p_1|}
    \end{align*}
    and  the point  $(p_0,\beta_0(p_0))$   is contained  in   the line
    defined  by  $(\beta_1(p_1),p_1)$ and $(i_0,i_1)$.   So the points
    $V(F_0\times  F_1)$, $(p_0,\beta_0(p_0))$ and $(\beta_1(p_1),p_1)$
    lie on a common hyperplane $H$.
  
    The product  $P_0\times P_1$ lies entirely  on one side of  $H$ by
    construction.  Suppose there  is a point $x$  of $S$ on the  other
    side of    $H$.  As $H$   is   a valid hyperplane   for  the ridge
    $F_0\times  F_1$, any  point beyond it  is  also beyond either the
    facet hyperplane of $F_0\times P_1$   or $P_0\times F_1$.   Assume
    the first.  For any $z\in S$ we have  either $\pi_0(z)\in S_1$, or
    $\pi_0(z)\in  V(P_0)$, or $\pi_0(z)\in V(E(P_0))\setminus V(P_0)$.
    $x\in  S$ is beyond  $F_0\times P_1$,  therefore only $\pi_0(x)\in
    V(E(P_0))\setminus V(P_0)$ is  possible.  Thus,  $\pi_0(x)$ is the
    unique  vertex of   $E(P_0)$ beyond $F_0$,  so  $\pi_0(x)=p_0$ and
    $x\in H$.

  \end{compactenum}
  \noindent
  \begin{minipage}{\textwidth}
    This  proves    that the two   conditions   (\ref{conditionA}) and
    (\ref{conditionB}) are sufficient for the existence of a polytopal
    realisation of $E(P_0\times P_1)$.\qedhere
  \end{minipage}
\end{proof}
In Section \ref{sec:realisations} we present some general applications
of Theorem  \ref{thm:construction}.    However, we mostly  use  a more
restrictive version      of   it.     We  tighten    the    conditions
(\ref{conditionA}) and (\ref{conditionB}) in the following way to make
them more manageable:

\noindent
\fbox{
  \parbox{\textwidth-.5cm}{
    \begin{compactenum}[(A$'$)]
    \item\label{conditionAp}      There      exist   vertex-preserving
      realisations of $E(P_0)$ and $E(P_1)$.
    \item\label{conditionBp} For $i=0,1$  let $T_i:=V(E(P_i))\setminus
      V(P_i)$.   There  are points  $s_i\in\interior(P_{i})$ and  some
      $0<r<1$ such that for any $v_0\in T_0$ and $v_1\in T_1$
      \begin{align*}
        r|s_0,v_0|=|s_0,q_0| \qquad\qquad (1-r)|s_1,v_1|=|s_1,q_1|,
      \end{align*}
      where $q_i$ is the  intersection  of the  segment from $s_i$  to
      $v_i$ and $|a,b|$ denotes the length  of the segment from $a$ to
      $b$ (Hence $\beta_i(x)\equiv s_i$ for $i=0,1$).
    \end{compactenum}
  }
}

\medskip
Theorem~\ref{thm:construction} immediately implies the following.
\begin{cor}\label{cor:specialcase}
  Let $P_0$, $P_1$ be a pair of polytopes with $\dim(P_0\times P_1)\ge
  3$   that    satisfies     {\rm  (\ref{conditionAp})}     and   {\rm
    (\ref{conditionBp})}.  Let $S\subset\R^{d_0}\times\R^{d_1}$ be the
  point set containing the following points:
  \begin{compactenum}[\rm (a)]
  \item all pairs $(p_0,p_1)$ for $p_0\in V(P_0)$, $p_1\in V(P_1)$,
  \item all pairs $(v_0,s_1)$ for $v_0\in T_0$,
  \item all pairs $(s_0,v_1)$ for $v_1\in T_1$.
  \end{compactenum}
  Then  $\conv(S)$  is  a vertex-preserving polytopal   realisation of
  $E(P_0\times P_1)$.\qed
\end{cor}
In this  setting, the only  connection between the construction of the
two factors is  the value $r$  of the ratio.  Thus,  if we construct a
polytope $P$  together with a  vertex-preserving realisation of $E(P)$
and a single point $s$ in its interior such that all segments from $s$
to  the  vertices of  $E(P)$ not in  $P$ intersect  $\boundary P$ with
ratio $r$, then we can combine this with any other such instance for a
ratio of $1-r$.

\section{Explicit realisations}

Now we apply the construction    of the previous section and   produce
examples  of   products  of polytopes  with  a  realisation   of their
$E$-polytope.   The  main   focus   is on  the  realisation   of   the
$E$-polytope $E_{mn}$ of a product  of an $m$-gon  and an $n$-gon.  We
produce polytopal realisations for all of them.  Afterwards we briefly
discuss examples in   dimensions $d\ge5$.  For  some of  the  examples
explicit data in the polymake format (Gawrilow and Joswig~\cite{GJ00})
are available on request.

\subsection{Products of polygons}\label{sec:polytopality} 

We  begin with products  of polygons and present a  method to obtain a
``flexible''  geometric realisation of  $E_{mn}:=E(C_m\times C_n)$ for
all $m,n\ge     3$.  We will   discuss degrees    of   freedom in this
construction in Section \ref{sec:realisations}.
\begin{theorem}\label{thm:EmnReal}
The CW spheres $E_{mn}$ are polytopal for all $m,n\ge 3$.
\end{theorem}
In the five  cases   when $m,n$ satisfy  $\frac{1}{m}+  \frac{1}{n}\ge
\frac{1}{2}$, polytopality   follows from a   construction of   Santos
\cite{Santos02},  \cite[Rem.~13]{Santos00}.   These  realisations  are
presented in Theorem \ref{thm:symmetry}.

\begin{figure}[b]
  \centering
  \begin{minipage}[b]{.33\textwidth}
    \centering
    \psfrag{E}{$E(\Delta)$}
    \psfrag{D}{$\Delta$}
    \psfrag{S}{$s$}
    \includegraphics[width=.99\textwidth]{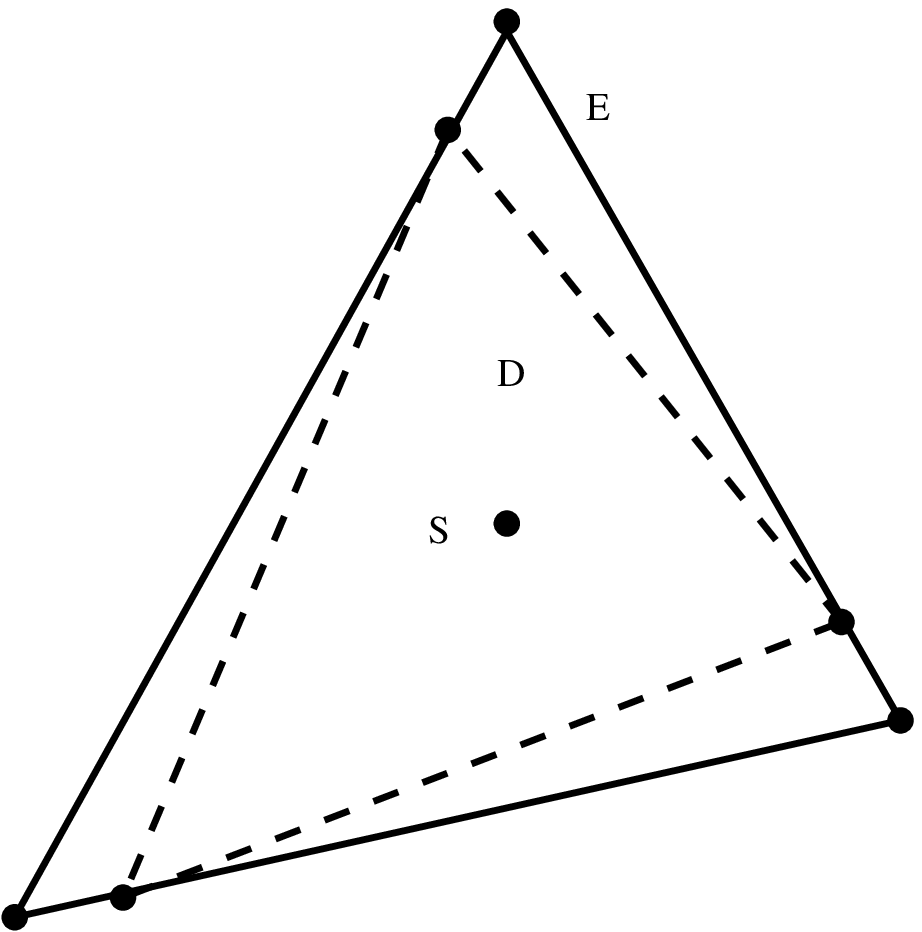}
    \caption{An example of a realisation of $D(3,\frac{1}{3})$.}
    \label{fig:dnr}
  \end{minipage}
  \hspace{.02\textwidth}
  \begin{minipage}[b]{.55\textwidth}
    \centering
    \psfrag{s}[b][b]{$s$}
    \psfrag{v0}[t][t]{$v_0$}
    \psfrag{v1m}[tr][tr]{$v_1^-$}
    \psfrag{v1p}[tl][tl]{$v_1^+$}
    \psfrag{v2m}[tr][tr]{$v_2^-$}
    \psfrag{v2p}[tl][tl]{$v_2^+$}
    \psfrag{v3}[b][b]{$v_3$}
    \psfrag{C}[br][br]{$C_6$}
    \psfrag{E}[b][b]{$E(C_6)$}
    \psfrag{pbv}[l][l]{$\overline p((v_2^+)_x)$}
    \psfrag{qbv}[r][r]{$\overline q((v_2^+)_x)$}
    \psfrag{pv}[l][l]{$p((v_1^+)_x)$}
    \includegraphics[width=.95\textwidth]{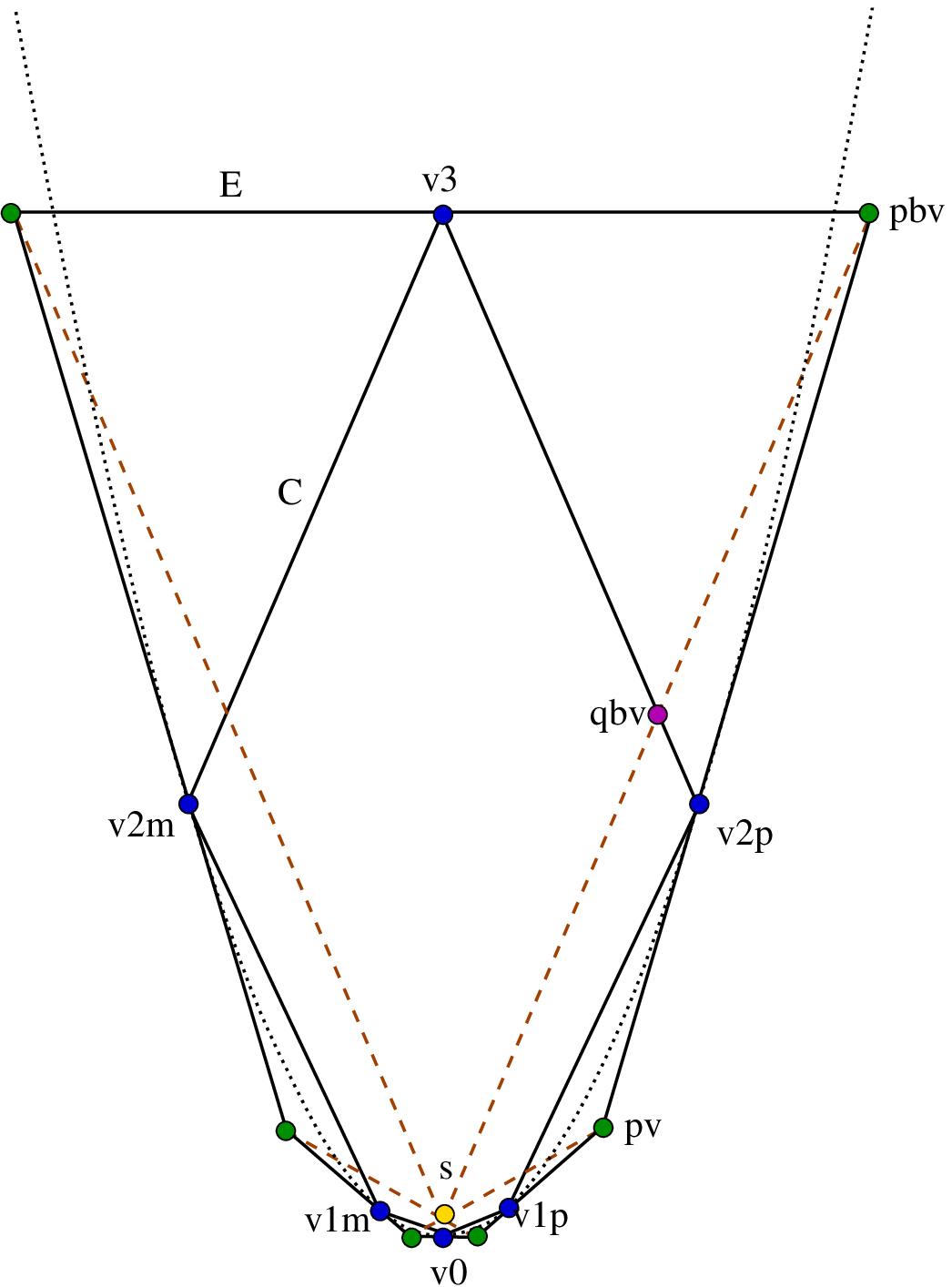}
    \caption{The  construction    of   $D(6,\frac35)$. Note  that  the
      $x$-axis has larger scale than the $y$-axis.}\label{fig:c-dmr}
\end{minipage}
\end{figure}

We use the  restricted setting of Corollary  \ref{cor:specialcase} for
the  proof and construct  only one of the two  factors, but subject to
the  conditions (\ref{conditionAp}$'$) and (\ref{conditionBp}$'$).  We
make the following definition for this.
\begin{defi}
  Let $\frac13<r<\frac23$. By  $D(k,r)$ we denote  a  realisation of a
  $k$-gon $C_k$ together with
\begin{compactitem}
\item a distinguished inner point $s$,
\item a  vertex-preserving $E$-polytope $E(C_k)$,  such that  segments
  from $s$ to vertices of $E(C_k)$  are intersected by the boundary of
  $C_k$ with ratio $r$.
\end{compactitem}
\end{defi}
See Figure  \ref{fig:dnr}  for an example.    To  realise $E_{mn}$  we
choose  a ratio $r$  between  $\frac13$ and $\frac23$  and combine the
points     of    $D(m,r)$    and       $D(n,1-r)$      according    to
Corollary~\ref{cor:specialcase}.   We  restrict  to $m\ge  4$   in the
following and refer to the proof of Theorem~\ref{thm:symmetry} for the
case $m=3$.

The next construction is   illustrated in Figure~\ref{fig:c-dmr}.  Let
$\Gamma$ denote the graph of the  parabola $x\mapsto x^2$ in the plane
$\R^2$  with coordinates $x$ and $y$.   We construct the polygon $C_m$
such that  all but one of the  vertices lie on $\Gamma$,  and $E(C_m)$
such that all but one of the edges are tangent to $\Gamma$.

\newcommand{\cint}{\mathcal C}
\newcommand{\cend}{\overline{\mathcal C}}
\newcommand{\eend}{\mathcal E}

Let $s$ be the point $(0,1)$ and define the three functions
\begin{align*}
  \cint(x)&:=\frac{x+\sqrt{(1-r)(2r+rx^2+x^2)}}{r}\\
  \cend(x)&:=\frac{x(x+\sqrt{x^2+2r^2x^2-2rx^2-2r+2r^2}-rx)}{r}\\
  \eend(x)&:=\frac{x+\sqrt{x^2+2r^2x^2-2rx^2-2r+2r^2}}{2r}.
\end{align*}
For any  $a\ge 0$ let $p(a)\in\R^2$  be the  intersection point of the
tangents  to $\Gamma$  in $(a,a^2)$  and  $(\cint(a),\cint(a)^2)$.  We
have the following facts about these functions.
\begin{lemma}\label{lemmafunctions}
  Let $0<r<\frac23$.
  \begin{compactenum}
  \item \label{l1} For any $a\ge 0$  the secant line between $(a,a^2)$
    and $(\cint(a),\cint(a)^2)$ intersects the segment between $s$ and
    $p(a)$ in a point $q(a)$ satisfying
    \begin{align*}
      r|s,p(a)|=|s,q(a)|,
    \end{align*}
    where $|x_0,x_1|$ denotes the length of the segment between $x_0$
    and $x_1$.
  \item  \label{l1a} For any $a\ge  1$ the line  between $(a,a^2)$ and
    $(0,\cend(a))$ intersects     the    segment   between     $s$ and
    $\overline{p}(a):=(\cend(a),\eend(a))$      in      a        point
    $\overline{q}(a)$ satisfying
    \begin{align*}
      r|s,\overline{p}(a)|=|s,\overline{q}(a)|.
    \end{align*}
  \item\label{l2} For any $a\ge 0$ we have $\cint(a)>1$.
  \item\label{l3}For any $a>1$ we have $\cend(a)>a^2$.
  \item\label{l4}       For    any    $a>1$    we     have  $\cend(x)=
    2x\eend(x)-x^2$.\hfill\qed
\end{compactenum}
\end{lemma}

With this  information at  hand   we are  ready to give   an iterative
construction for $D(m,r)$ in the case  $m\ge 4$ and $0<r<\frac23$.  We
distinguish the two cases $m$ even and $m$ odd.
\begin{compactitem}
\item For $m$  even we choose  the point  $v_0^+:=(0,0)$ as  our first
  vertex  of    $C_m$,   for   $m$ odd   we     take the   two  points
  $v_0^\pm:=(\pm\sqrt{\frac{1-r}{1+r}},\frac{1-r}{1+r})$.
\item   In     the  $i$-th   step   we   extend    with   the   points
  $v_i^\pm:=(\pm\cint(a),\cint(a)^2)$, where $a$ is the $x$-coordinate
  of $v_{i-1}^+$.
\item  We  repeat the previous step   until we  have constructed $m-1$
  points of $C_m$.
\item   In    the       last     step    we   add       the     vertex
  $v_{\lfloor\frac{m}{2}\rfloor}:=(0,\cend(a))$,  where $a$ is     the
  $x$-coordinate of $v_{\lfloor\frac{m}{2}\rfloor-1}$.
\end{compactitem}
The edges of $E(C_m)$ are the tangents to $\Gamma$  in the vertices of
$C_m$, except  for $v_{\lfloor\frac{m}{2}\rfloor}$,  where we take  the
horizontal line  running through $v_{\lfloor\frac{m}{2}\rfloor}$. This
line intersects the tangents to $(v_{m-1}^\pm,(v_{m-1}^\pm)^2)$ in the
points $(\pm \eend(v_{m-1}^+),\cend(v_{m-1}^+))$.

Lemma~\ref{lemmafunctions}   guarantees      the  condition on     the
intersection ratio.   The   point   $s=(0,1)$  is inside     $C_m$  by
Lemma~\ref{lemmafunctions}(\ref{l2}),(\ref{l4}).     Hence      we can
construct   $D(m,r)$ for  $0<r<\frac23$  and  $m\ge   4$.  Finally,  a
realisation for  $m=3$  and $\frac14<r<1$  is  given  in the proof  of
Theorem \ref{thm:symmetry}.  By combining $D(n,r)$ with $D(m,1-r)$ for
some    $m\ge 3$  we   obtain     $E_{mn}$.    This  proves    Theorem
\ref{thm:EmnReal}.

\begin{remark}
  $D(m,r)$ can in fact be constructed for any  $r$ between $0$ and $1$
  and any $m\ge 3$,  but the formulas for  the vertices and the  cases
  one has  to distinguish in the construction  tend to get complicated
  rather quickly.
\end{remark}

\subsection{Some higher-dimensional examples}

Satisfying the conditions (\ref{conditionA}) and (\ref{conditionB}) is
more difficult  if the two factors $P_0$  and  $P_1$ have more facets.
Thus  in higher dimensions  and for ``more  complex'' polytopes, it is
usually hard to  find maps $\beta_0$ resp.\  $\beta_1$, unless one can
exploit some kind of symmetry.

There are,  however, two  obvious  families of polytopes  that we  can
choose as factors  of  a product polytope: The  $d$-simplex $\Delta_d$
and   the $d$-cube $C_d$.   Both can  be  realised together with their
$E$-construction  satisfying   even  the restrictive    conditions  of
Corollary \ref{cor:specialcase}.

\newcommand{\simplex}{\Delta}

\begin{figure}[b]
    \centering
    \psfrag{v}{$v$}
    \psfrag{H}[t][t]{$H$}
    \psfrag{Dp}[r][r]{$\Delta'$}
    \psfrag{D}[r][r]{$\Delta$}
    \psfrag{F}[r][r]{$F$}
    \psfrag{Fp}[tl][tl]{$F'$}
    \includegraphics[height=.33\textheight]{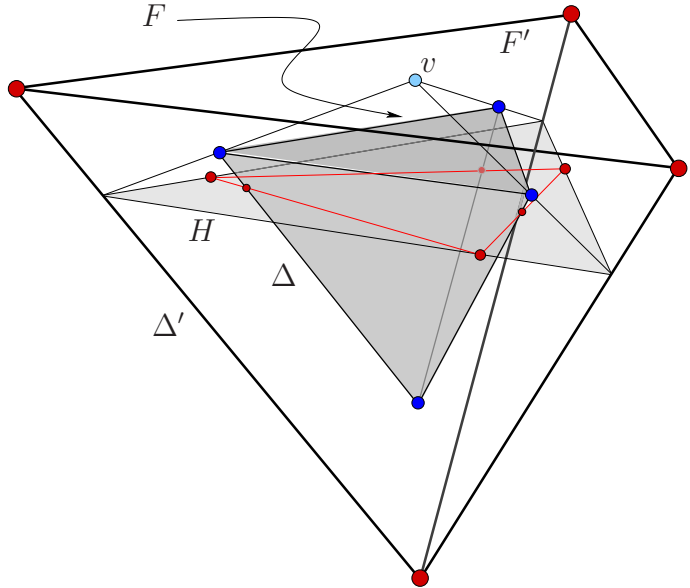}
    \caption{The construction for    a $3$-simplex. The   solution for
      $d=2$ used in the plane $H$ is indicated with thin lines.}
    \label{fig:construct3simplex}
  \end{figure}

\begin{compactitem} 
\item The cube with its $E$-construction  and an intersection ratio of
  $r=\frac12$ can  be  realised  as follows:   For  $C_d$ we  take the
  standard $\pm1$-cube. The new vertices for the $E$-polytope are $\pm
  2\cdot e_i$,  where $e_i$ are the standard  unit basis  vectors. The
  origin is an inner point $s$ satisfying all requirements.

\item  The construction  for the $d$-simplex  $\simplex_d$ is slightly
  more difficult.    We give an  inductive construction  that produces
  realisations  for any ratio  $\frac{1}{2}\le r<1$, that is, at least
  half of the segment is inside $\simplex$ (where $r$ is the parameter
  appearing      in      the   conditions    in    the       box    on
  page~\pageref{prod-suff-cond}).   We can clearly  construct  such  a
  realisation for a triangle, i.e.\ for a simplex of dimension $d=2$.
   
  For $d>2$ we  take a regular  realisation $\simplex$ of the  simplex
  and a scaled  version $\simplex':=\frac{1}{r}\cdot\simplex$ with the
  same barycentre.  We  choose one  facet  $F$ of  $\simplex$ and  the
  corresponding scaled facet $F'$ in $\simplex'$.  Place the first new
  vertex $v$ in the barycentre of $F'$.  The vertices of any ridge $R$
  of $F$ together with the point $v$ uniquely define a hyperplane. $F$
  has  $d$    ridges,  so   we obtain     $d$   different  hyperplanes
  $H_1,\ldots,H_d$  by this.   $H_1,\ldots,H_d$   intersect  all facet
  hyperplanes  of    $\simplex'$,      except   that  to   $F'$,    in
  co\-di\-men\-sion-2-planes that lie in a  common hyperplane $H$. $H$
  is parallel to $F$.
  
  $H$  cuts    $\simplex$     and $\simplex'$    in    two   simplices
  $\tilde\simplex$ and  $\tilde\simplex'$ of dimension $d-1$.  (Recall
  that  $r\ge\frac{1}{2}$,   so $H$  intersects $\simplex$   below the
  barycentre if viewed from $F$.)  $\tilde\simplex'$ is (viewed in the
  hyperplane $H$) a scaled version of  $\tilde\simplex$ with a scaling
  factor  $\frac{1}{r'}\le   \frac{1}{r}$.  By   induction,  we have a
  solution for  the  corresponding  problem for  $\tilde\simplex$  and
  $r'\ge r\ge   \frac{1}{2}$  in $H$  (where the  inner  point  is the
  projection of the barycentre of $\simplex$).
  
  These points, together with the one vertex $v$ chosen before, give a
  realisation of  $E(\simplex)$  that   satisfies the   conditions  of
  Corollary~\ref{cor:specialcase}.                                 See
  Figure~\ref{fig:construct3simplex} for the case $d=3$.
\end{compactitem}

We can combine  such a simplex or cube  with another simplex, cube  or
some $D(n,\frac{1}{2})$ to obtain the $E$-polytope of this product.

\section{Properties of the family $E_{mn}$}

This section collects several properties of the polytopes $E_{mn}$. In
particular we count degrees of freedom for the realisation of $E_{33}$
and   prove  that not all  combinatorial    symmetries of $E_{mn}$ are
geometrically realisable.

\subsection{Self-duality}

The   polytopes $C_m\times  C_n$   are  simple,   thus   we know  from
\cite[Thm.\ 1.6]{PZ03} that  $E_{mn}$ is $2$-simple and $2$-simplicial
(that is,  all  $2$-faces are   triangles and  all  edges  are  in $3$
facets).  In  particular the $f$-vector of  $E_{mn}$ is symmetric (cf. 
Eq.  (\ref{eq:fvector})).

The polytopes $E_{mn}$  are in fact  self-dual.  This is  not true for
arbitrary $2$-simple, $2$-simplicial polytopes, which can be seen from
the hypersimplex $E(\Delta)$ obtained  from the  $4$-simplex $\Delta$. 
This polytope has a  facet-transitive automorphism group acting on its
10 bipyramidal   facets, while   the  dual  has   5 tetrahedral and  5
octahedral facets.

\begin{theorem}[Ziegler \cite{Ziegler03}]
  Each of the polytopes $E_{mn}$ ($n,m\ge 3$) is self-dual, with an
  anti-automorphism of order $2$.
\end{theorem}

\begin{proof}
  Number   the vertices   of    an  $k$-gon  $C_k$  consecutively   by
  $v_0,\ldots,v_{k-1}$. We  take indices modulo  $k$.  The vertices of
  $C_m\times C_n$ are   $v_{i,j}:=(v_i,v_j)$ for $0\le  i\le m-1$  and
  $0\le j\le n-1$.  We have two types of facets in the product:
\begin{align*}
  F'_i&=\conv(\{v_{ij},v_{i+1,j}\mid j=0,\ldots, n-1\})\\
  F''_j&=\conv(\{v_{ij},v_{i,j+1}\mid i=0,\ldots, m-1\}).
  \intertext{We denote the new vertex beyond  $F'_i$ by $v'_i$ and the
    one beyond  $F''_i$ by $v''_i$. The  facets  of $E(C_m\times C_n)$
    are now of the form}
  G_{ij}&=\conv(v_{ij},v_{i+1,j},v_{i,j+1},v_{i+1,j+1},v'_i,v''_j)\\
  G'_i&=\conv(\{v_{ij}\mid j=0,\ldots, n-1\},v'_{i-1},v'_i)\\
  G''_j&=\conv(\{v_{ij}\mid     i=0,\ldots,    m-1\},v''_{j-1},v''_j).
  \intertext{From    this we can  read  off   the facets   a vertex is
    contained in:}
  v_{ij}&\in G_{ij},G_{i-1,j},G_{i,j-1},G_{i-1,j-1},G'_{i-1},G''_{j-1}\\
  v'_i&\in G'_{i},G'_{i+1}, G_{ij}\text{ for } j=0,\ldots, n-1\\
  v''_j&\in G''_{j},G''_{j+1}, G_{ij}\text{ for } i=0,\ldots, m-1.
\end{align*}
Thus the following correspondence gives a self-duality of order $2$: 
\begin{align*}
G_{ij}&\longleftrightarrow v_{-i,-j}\qquad
G'_{i} \longleftrightarrow v'_{-i}\qquad
G''_{j}\longleftrightarrow v''_{-j}\qedhere
\end{align*}
\end{proof}

For    $m=n$, this  result was    obtained    previously by  G{\'e}vay
\cite{Gevay04}.

\begin{remark}
  There are examples of $3$-polytopes that are self-dual, but that do not
  have a self-duality of order $2$ (cf.\ \cite[Ex.~3.4.3,
  p.52d]{Gruenbaum67}).
\end{remark}

\subsection{$E_{mn}$ constructed from regular polygons}\label{sec:symmetry}

Only  in a few cases there  are ``more symmetric'' realisations of the
polytopes $E_{mn}$: We prove that there are only five choices of pairs
$(m,n)$ (up  to interchanging   $m$  and $n$) such  that  we  can take
regular polygons as  input for the  construction described  in Theorem
\ref{thm:construction}  in the    restricted version   of    Corollary
\ref{cor:specialcase}. We will see in the next section that these five
cases   are also the only   cases in which  the  product of two cyclic
groups induced from rotation of the vertices in the two factors can be
a subgroup of the geometric symmetry group.  The next theorem is based
on Santos~\cite[Rem.~13]{Santos00}, \cite{Santos02}.

\begin{theorem}\label{thm:symmetry}
  There  are polytopal realisations of  $E_{mn}$  for which projection
  onto the first and last two coordinates yields in both cases
  \begin{compactenum}[\rm (1)]
  \item  regular polygons for the  polygon in  $C_m\times C_n$ and its
    $E$-con\-struc\-tion,
  \item and all intersection ratios are equal in each factor
    \label{ass:ratioequal}
  \end{compactenum}
  if and only if $\frac{1}{m}+\frac{1}{n}\ge\frac{1}{2}$.
\end{theorem}

\begin{figure}[b]
\centering
\begin{minipage}{.8\textwidth}
\psfrag{s}{$l$}
\psfrag{P}{$P$}
\psfrag{EP}{$E(P)$}
\centering
\includegraphics[width=3.3cm]{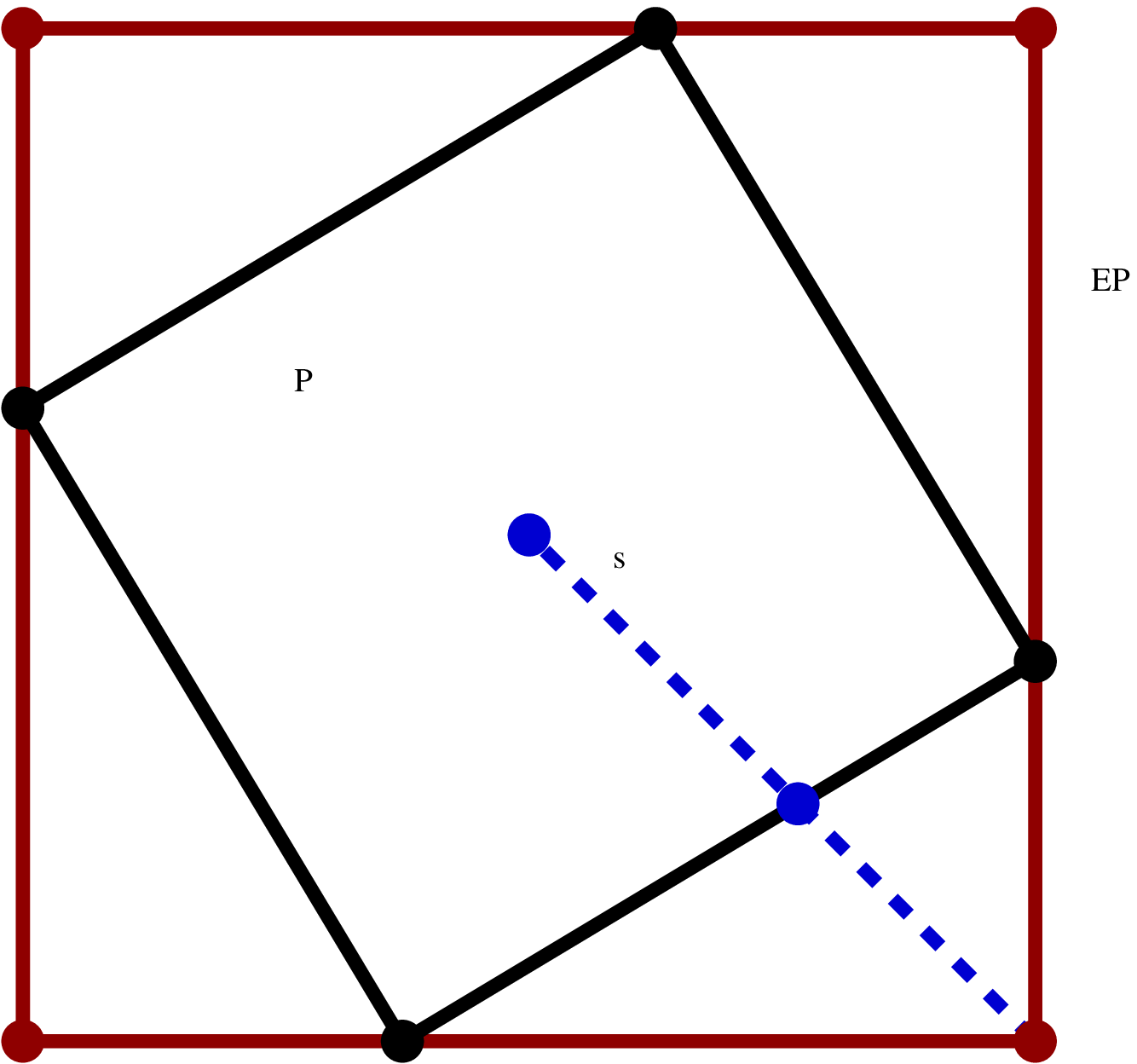}
\hspace{2cm}
\includegraphics[width=3cm]{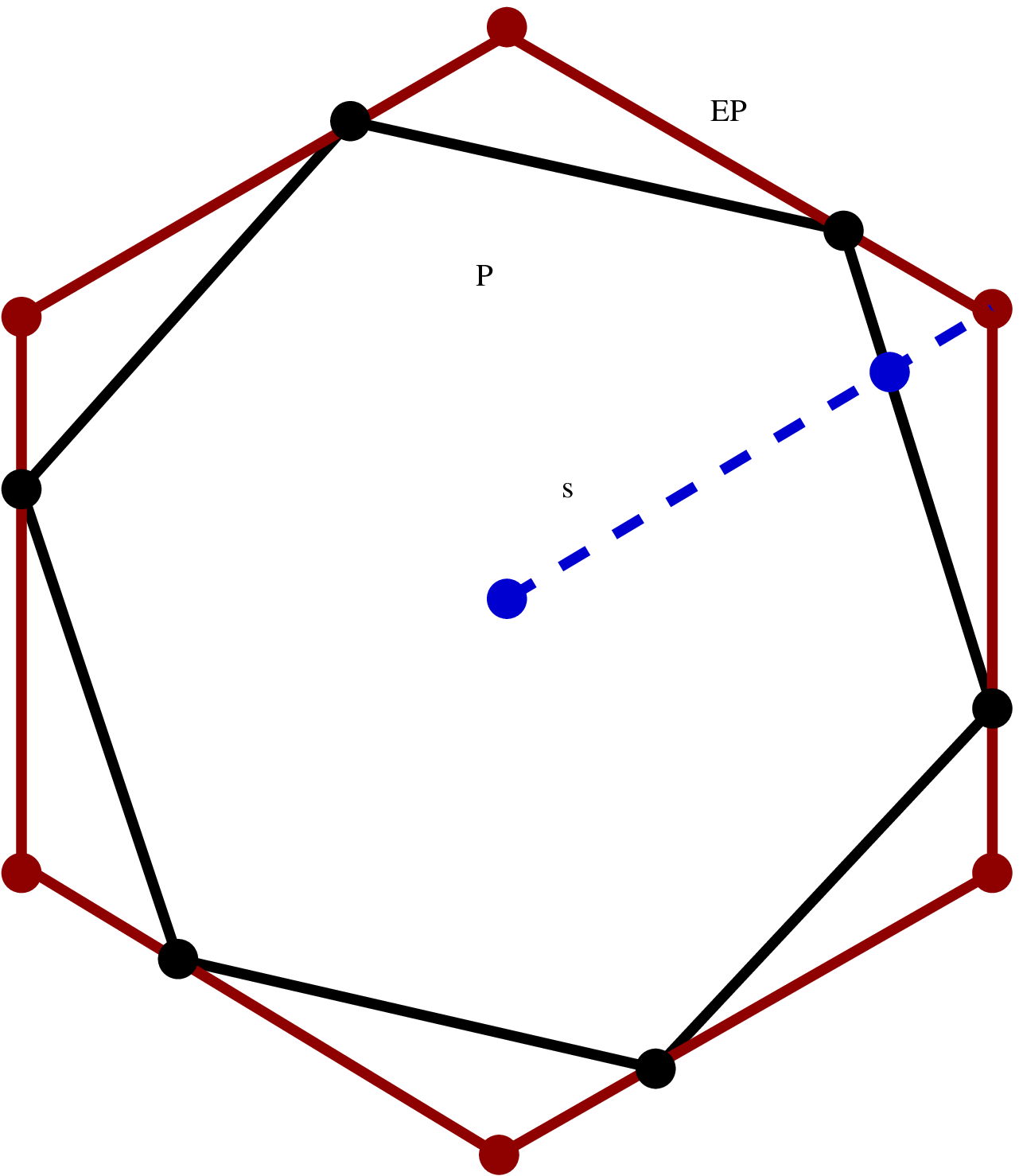}
\caption{Two projections that 
  satisfy the restrictions of Theorem~\ref{thm:symmetry}}
\label{fig:SymmProj}
\end{minipage}
\end{figure}
\begin{proof}
  The condition  on  the ratio implies  that  the images of   the maps
  $\beta_0$ and $\beta_1$  appearing  in the construction  of $E_{mn}$
  are single points in  the interior of the  polygons $C_m$ and $C_n$.
  These points must  be the barycentres  if the polygons are  regular.
  We may assume that this is the origin.
  
  We can now  generate all configurations  of a  regular polygon $C_m$
  together with $E(C_m)$  in the following  way: Start with  a regular
  polygon $E(C_m)$ centred at the origin and choose a vertex for $C_m$
  on each of the  edges.  As $C_m$ is  regular, the vertices of  $C_m$
  divide  each edge with  equal  ratio.   The  segments  considered in
  (\ref{conditionB}) are  the  segments $l$  between  the origin and a
  vertex of $E(C_m)$.  We are interested in the possible values of the
  ratio with which they are intersected by the edges of $C_m$.
  
  Choosing  the vertices of $C_m$  close to  those  of $E(C_m)$ we see
  that we can have an arbitrarily high portion of $l$ inside $C_m$. On
  the other hand, the portion inside  $C_m$ is minimised when we place
  the vertices of $C_m$ in the centre of the edges.  In this case, the
  fraction   of $l$  outside  $C_m$  is  $\sin^2(\frac{\pi}{n})$.   By
  condition  (\ref{conditionB}),  the fraction  lying outside  for one
  polygon  and its $E$-construction   has to match the  fraction lying
  inside for the other polygon. This gives the following inequalities:
  \begin{align*}
    1-\sin^2\left(\frac{\pi}{m}\right)\le\sin^2\left(\frac{\pi}{n}\right)
    \qquad\text{and}\qquad
    1-\sin^2\left(\frac{\pi}{n}\right)\le\sin^2\left(\frac{\pi}{m}\right),
  \end{align*} 
  which are equivalent to the condition given in the theorem.
\end{proof}
We can  determine  all possible  values  for  the inequalities of  the
theorem explicitly.

\begin{cor}\label{cor:symmetry}
  There are realisations of  $E_{mn}$ from regular polytopes only  for
  the following pairs $(m,n)$ (up to interchanging $m$ and $n$):
  \begin{align}
    (3,3)\;,\;(3,4)\;,\;(3,5)\;,\;(3,6)\;,\;(4,4)\tag*{\qed}
  \end{align}
\end{cor}

\begin{remark}
  We    made    assumption    (\ref{ass:ratioequal})    in     Theorem
  \ref{thm:symmetry} mainly  because this is  the case we need  in the
  next section.  A less restrictive version of ``symmetry'' would only
  require the points in the images  of $\beta_0$ and $\beta_1$ to also
  form regular polygons (if we take the points in the order induced by
  the $E$-construction of the other factor).  For small $m=n$ this has
  solutions where  all   points  in the images   are  different.   See
  Table~\ref{tab:less_symmetry} for an example  of an  $E_{44}$.  Note
  however,  that   this  severely  reduces  the   number  of geometric
  symmetries compared to the case of the theorem.
\end{remark}

\subsection{Combinatorial versus geometric symmetries}

There are two different notions of symmetry for a polytope $P$.
\begin{defi}
  Let $P$ be a polytope with a given geometric realisation. Any affine
  transformation  $T$ of the ambient  space that preserves $P$ set-wise
  is called  a {\em geometric symmetry  transformation}.  The group of
  all  such  transformations  is called  the  {\em  geometric symmetry
    group}.
  
  To any polytope $P$ we  can associate the poset of  all faces of $P$
  ordered by inclusion.  This is called the  {\em face lattice} $F(P)$
  of $P$.  A {\em combinatorial symmetry} of $P$ is an automorphism of
  $F(P)$.  The  group  of all  combinatorial  symmetries  is  the {\em
    combinatorial symmetry group}.
\end{defi}
The  combinatorial symmetry  group  is  independent of  a realisation,
while the geometric symmetry group highly depends on the choice of the
realisation.

A geometric symmetry    maps  $k$-faces to  $k$-faces   and  preserves
incidences.  Therefore any  geometric symmetry induces a combinatorial
symmetry.  On the other hand, a combinatorial symmetry in general does
not induce a geometric one. However, there are not many examples known
of polytopes where these two groups  differ for all possible geometric
realisations of a  polytope.   Bokowski, Ewald, and Kleinschmidt  have
provided a $4$-dimensional  example on $10$  vertices in \cite{BEK84}.
Dimension $4$ is smallest possible for such examples,  as it is known,
that for $3$-polytopes, and  for polytopes  with  few vertices  in any
dimension,    there  are    realisations for     which   geometric and
combinatorial symmetry group coincide (see \cite{Mani71} for the first
and  \cite[p.120]{Gruenbaum67} for the   second result).  We show that
our product construction provides an  infinite series of $4$-polytopes
with non-realisable  geometric symmetries.   We construct an  explicit
example of such a symmetry for  the proof.  Previously it was observed
by G{\'e}vay that there is no  polytopal realisation of the CW spheres
$E_{mm}$ with the full symmetry group, except in the case $m=4$.  This
is also a consequence of Corollary \ref{cor:norot} below.

\begin{theorem}\label{thm:nosym}
  For relatively  prime  $m,n\ge5$  all $E_{mn}$  admit  combinatorial
  symmetries   that  cannot  be    realised  as    affine  symmetry
  transformations of a geometric realisation of $E_{mn}$.
\end{theorem}

\begin{table}[bt]
\setlength{\fboxsep}{5pt}
\framebox{\parbox{\textwidth-2\fboxsep}{{\sf 
{\bfseries The three symmetries involved in the proof of Theorem \ref{thm:nosym}}

\medskip
Notation:
\begin{compactitem}
\item $v_0, v_1, v_2$: vertices of $C_3$ 
\item $w_0, w_1, w_2, w_3$ vertices of $C_4$
\item $e(j)$: edge from vertex number $j$ to $j+1$ (mod $3$ or $4$) in both polygons. 
\end{compactitem}

\medskip
Number the vertices $p_k$ of $P_{34}$ in the following way:
\begin{list}{}{\settowidth{\labelwidth}{\text{$15\le k\le 19$:}}\advance\leftmargin\labelwidth\setlength{\itemsep}{-3pt}}
\item[$ 0\le k\le 11$:] vertices $(v_{k\text{ div } 4}, w_{k-4(k\text{
      div } 4)})$
\item[$12\le k\le 14$:] vertices added above $e(k-12)\times C_4$
\item[$15\le k\le 19$:] vertices added above $C_3\times e(k-15)$
\end{list}
 Then the combinatorial symmetries are given as (permutation notation, vertex numbers of $p_k$):
\begin{align*}
\tilde S_3&:=(0, 4, 8)(1, 5, 9)(2, 6, 10)(3, 7, 11)(12, 13, 14)(15)(16)(17)(18)\\
\tilde S_4&:=(0, 1, 2, 3)(4, 5, 6, 7)(8, 9, 10, 11)(12)(13)(14)(15, 16, 17, 18)\\
T&:=(0, 5, 10, 3, 4, 9, 2, 7, 8, 1, 6, 11)(12, 13, 14)(15, 16, 17, 18)
\end{align*}}}}
\caption{The combinatorial symmetries $\tilde S_3$, $\tilde S_4$, and $T$ acting on $P_{34}$.}
\label{tab:symmetries}
\end{table}
\begin{proof}
Let   $P_{mn}$   be   any    geometric  realisation  of   a   polytope
combinatorially  equivalent  to an  $E_{mn}$.   Seen  as  a PL sphere,
$P_{mn}$ can still  be viewed as the   result of the  $E$-construction
applied   to a  PL  sphere  which  is combinatorially  equivalent to a
product of two polygons.

Here is a non-realisable combinatorial  symmetry $T$ of $P_{mn}$.  Let
$C_m$  and $C_n$ denote  polygons  with vertices $v_0,\ldots, v_{m-1}$
resp.\  $w_0,\ldots, w_{n-1}$   numbered  in  cyclic order.  We   take
indices modulo $m$ resp.\ $n$.  Let $S$  be the combinatorial symmetry
of a polygon that maps the $j$-th to the $(j+1)$-th vertex.

$S$ induces   a combinatorial symmetry   $S_m$  on $C_m\times  C_n$ by
mapping  a  vertex $(v_i,w_j)$ to  $(v_{i+1},w_j)$  for any $0\le j\le
m-1$.  Similarly $S$ induces a symmetry $S_n$ shifting the vertices of
$C_n$.  Both   symmetries uniquely extend  to combinatorial symmetries
$\tilde S_m$ and $\tilde  S_n$ of $E(C_m\times  C_n)$.  Let $T$ be the
combinatorial symmetry of $P_{mn}$ obtained  by first applying $\tilde
S_m$ and then  $\tilde  S_n$.  See Table~\ref{tab:symmetries}  for  an
example.

A geometric  realisation   of $P_{mn}$   need  not  have  the  product
structure      induced   by     the     construction    of     Theorem
\ref{thm:construction}.  However,  by  looking at vertex  degrees, for
$m,n\ge 5$ we can decide which  of the vertices of $P_{mn}$ ``belong''
to the product  and which  are ``added''  by the  $E$-construction:  A
vertex  of the product  always has  degree $8$,  as $C_m\times C_n$ is
simple, so any vertex has four neighbours and is in  four facets.  The
added vertices all have degree $2m$ or $2n$. Denote the vertex sets by
$V_p$ and $V_e$.

The  proof is  roughly as   follows.   Suppose there  is  a  geometric
realisation $T_g$ of $T$  for some $P_{mn}$.   First we prove that any
$P_{mn}$ having $T_g$  as  a geometric symmetry has  the  form of  the
construction  in   Theorem \ref{thm:construction}.   Then the symmetry
implies   that both   factors are  of   the  form  defined  in Theorem
\ref{thm:symmetry}.  Corollary  \ref{cor:symmetry} finally   tells  us
that for $m,n\ge 5$ there are no such realisations.

As $T_g$ set-wise fixes the the vertices of $P_{mn}$ it also fixes the
centroid of the vertices of $P_{mn}$.  After a suitable translation we
can assume that $T_g$ is a linear  transformation.  As $m$ and $n$ are
relatively   prime, there is a   $k_m\in\N$ such that $T_m:=T_g^{k_m}$
restricted to the set $V_p$ acts as $\tilde S_m$. Similarly there is a
$k_n$ such that $T_n:=T_g^{k_n}$  reduces to a realisation  of $\tilde
S_n$.  $T_m$ and $T_n$ are linear transformations.

By construction $P_{mn}$ has   two  different combinatorial types   of
facets: Bi\-py\-ra\-mids over an $m$-gon and over an $n$-gon.  For any
facet we call the vertices of the polygon (i.e.\ those vertices of the
facet belonging to $V_p$) the base vertices.

Let $F$ be a facet of $P_{mn}$ of  the first type.  The symmetry $T_m$
shifts the base vertices by one and fixes the two apices.  Thus, $T_m$
also fixes  the centroid $c_F$ of the  base vertices  of $F$ and $T_m$
restricted to the  hyperplane $H_F$    defined  by $F$  is a    linear
transformation $T_m^F$ in   $H_F$ (if we  put  the origin of  $H_F$ in
$c_F$).   Now $T_m$ fixes  the two apices  of $F$  and  thus fixes the
whole  line through the  apices.  So $T_m^F$ splits  into a map fixing
the axis and a linear transformation  of a two dimensional transversal
subspace.   The axis must  contain $c_F$   and the  subspace  the base
vertices of $F$.  So the   base vertices of $F$   lie in a common  two
dimensional affine subspace  of $\R^4$.  Similarly, the  base vertices
of any other bipyramidal facet with a base equivalent  to $C_m$ lie in
a common $2$-plane.  These $2$-planes  are set-wise preserved by $T_m$
and therefore must be parallel.

The same argument proves   that  all bases of facets   combinatorially
being bipyramids  over $n$-gons do lie in  parallel $2$-planes.  These
$2$-planes  must  be  transversal  to  the  $2$-planes  containing the
$m$-gons: Otherwise  the    vertices in  $V_p$   all lie   in  a three
dimensional  subspace.  As $P_{mn}$ is four  dimensional, at least one
of the vertices of $V_e$ has to lie outside this $3$-space.  But there
are no edges between vertices in $V_e$.

Applying an  appropriate   linear transformation to  $P_{mn}$   we can
assume that the $2$-spaces containing the $m$-gons are parallel to the
$x_1$-$x_2$-plane and the  ones containing the  $C_n$ are  parallel to
the $x_3$-$x_4$-plane.   $T$ rotates the copies  of $C_m$ in $P_{mn}$,
so they must all be equivalent.  Similarly, all the polygons $C_n$ are
equivalent.      So  $P_{mn}$    is      an  instance    of    Theorem
\ref{thm:construction}.

Consider again  the facet $F$  with base equivalent  to  $C_m$ and the
restricted map  $T_m^F$.  Further restricting  $T_m^F$ to the subspace
containing  the base vertices  defines  a linear  map $T_b$  on $\R^2$
shifting  the vertices  of  a polygon.  So  $T_b$  generates a  finite
subgroup  of $Gl(2,\R)$ and therefore  must be conjugate to an element
of    $O(2,\R)$   (cf.      Schur~\cite{Schur11},      see        also
McMullen~\cite{McMullen68}).  The same argument applies to facets with
base  $C_n$.  As  the copies of  $C_m$  and  $C_n$ lie  in transversal
subspaces of $\R^4$, we can apply the  conjugation for $C_m$ and $C_n$
simultaneously and therefore both polygons are regular up to an affine
map.

Finally look at the $n$ vertices added above facets of $P_{mn}$ of the
type $C_m\times  F$ for  an  edge $F$ of  $C_n$.   Projecting onto the
$2$-space of $C_m$ they  lie inside $C_m$ (they  form the set $S_1$ in
the  construction of Theorem  \ref{thm:construction}).  They are fixed
by the symmetry  $\tilde S_m$.  As this map  has only one fixed point,
the points  in  $S_1$ must  coincide.  The same  applies  to the added
vertices above facets of type $F\times C_n$.   (Note that, even though
$T$ is a  symmetry  of the $E_{44}$  in Table~\ref{tab:less_symmetry},
the map $\tilde S_4$ is  not,  and cannot be   obtained as a power  of
$T$.)

Now  we are in the  situation described in Section \ref{sec:symmetry}. 
But according  to Corollary  \ref{cor:symmetry} this  can only be  the
case if  at least one of  $m$ and $n$ is  less than  $5$.  This proves
Theorem \ref{thm:nosym}.
\end{proof}

The  same   argument  also proves   that  Corollary \ref{cor:symmetry}
describes  all possible cases for which  $P_{mn}$ can have the product
$\Z_m\times \Z_n$ of two cyclic groups  induced by the rotation of the
vertices in the  two polygons as a  subgroup of its geometric symmetry
group.  In this case we do not need $m$ and $n$ to be relatively prime
as the   two  symmetries $\tilde S_m$   and  $\tilde  S_n$  itself are
contained in $\Z_m\times \Z_n$ acting on $P_{mn}$.
\begin{cor}\label{cor:norot}
  The combinatorial symmetry group of $E_{mn}$ contains a subgroup $G$
  isomorphic  to $\Z_m\times \Z_n$  induced  by  rotation  in the  two
  polygon factors.
  
  The geometric symmetry group of  a polytope $P_{mn}$ combinatorially
  equivalent  to $E_{mn}$ can contain  a subgroup  inducing $G$ on the
  face lattice only for $(m,n)= \{(3,3), (3,4), (3,5), (3,6), (4,4)\}$
  (up to interchanging $m$ and $n$).  \qed
\end{cor}
Hence, in particular, $E_{44}$ and $E_{36}$ are the only two polytopes
that have  a  geometric   realisation   realising all    combinatorial
symmetries.

\begin{remark}
  G{\'e}vay \cite{Gevay04} pointed out that along the lines of Theorem
  \ref{thm:nosym}    one can also   prove  that  the only  ``perfect''
  polytopes  among the  realisations of the   $E_{mn}$ are the regular
  $24$-cell and $E_{33}$ constructed  as in Theorem \ref{thm:symmetry}
  with intersection ratio $r=1/2$.   A rough definition of perfectness
  is as   follows:  A geometric realisation    $P$ of  a   polytope is
  \emph{perfect} if all  other  geometric realisations having,  up  to
  conjugation   with  an  isometry,  the same   subset   of the affine
  transformations as symmetry group,   are already isometric  to  $P$. 
  See \cite{Gevay02} for a precise definition.
\end{remark}

\subsection{Realisation spaces of $E_{33}$ and $E_{44}$}\label{sec:realisations}

We  determine the degrees  of  freedom that we have   in the choice of
coordinates for $E_{33}$. We consider two  realisations to be equal if
they  only differ  by a projective  transformation.  Thus, we  will be
interested in the dimension of the following spaces.
\begin{defi}
  The {\em realisation space} of  a $d$-polytope $P$ with $n$ vertices
  is  the space $\rs(P)$  of  all sets of   $n$ points in $R^d$  whose
  convex hull  is  combinatorially equivalent to   $P$. $\rs(P)$ is  a
  subset of $\R^{d\cdot n}$.
  
  The {\em projective realisation space} $\rs_{proj}(P)$ of a polytope
  is the space  of all possible  geometric realisations of a polytope,
  up to projective equivalence.  It  is the quotient space of $\rs(P)$
  where   two  realisations are  equivalent  if there  is a projective
  transformation mapping one onto the other.
\end{defi}
We work  out  the case of   $E_{33}$ explicitly and  present a simple
$4$-parameter   family of  $E_{44}s$.    We  prove  that  this  family
intersects  four  different  equivalence  classes  of the   projective
realisation space  $\rs_{proj}(E_{44})$.  Therefore, this  space is at
least four dimensional.

\subsubsection{The realisation space of $E_{33}$.}

The vertex  sets of all   realisations of $E_{33}$ obtained  from  the
construction in  Theorem \ref{thm:construction} contain the vertex set
of  an  orthogonal product  $C_3\times  C_3$ of   two triangles.  This
reduces the number   of  possible degrees of  freedom  compared  to an
arbitrary realisation. The  next  theorem determines the dimension  of
the  space  of all  realisations of   $E_{33}$  that are  projectively
equivalent  to a realisation  containing such an orthogonal product.

\begin{theorem}\label{thm:relspace}
  $\dim(\rs_{proj}(E_{33}))\ge 9$.
\end{theorem}

Before we   prove this   we introduce   a special  way  to   construct
realisations  of two triangles  and their $E$-polytopes satisfying the
conditions (\ref{conditionA}) and (\ref{conditionB}).

\begin{theorem}\label{thm:triangleratios}
  Given two (arbitrary)  triangles $\Delta$ and  $\Delta'$ there is an
  open subset $R$ in $\R^9$ such that, if we  take the nine entries of
  a   vector  in that  set  as  the  nine   ratios appearing  in  {\rm
    (\ref{conditionB})} (in some previously  fixed order), then  there
  is a realisation of $E_{33}$ having these intersection ratios.
\end{theorem}
\begin{proof}
  This is basically proven by  describing a realisation as a  solution
  of a set of linear equations, but we have to introduce some notation
  to write down the equations.
 
  \begin{figure}[b]
    \begin{minipage}{\textwidth}
      \centering
      \psfrag{a}{$s_a$}  
      \psfrag{b}{$s_b$}
      \psfrag{c}{$s_c$}
      \psfrag{ga}{$g_a$} 
      \psfrag{gb}{$g_b$}
      \psfrag{gc}{$g_c$}
      \psfrag{laap}{$l_{aa'}$}
      \psfrag{labp}{$l_{ab'}$}
      \psfrag{lacp}{$l_{ac'}$}
      \psfrag{lcap}{$l_{ca'}$}
      \psfrag{lcbp}{$l_{cc'}$}
      \psfrag{lccp}{$l_{cb'}$}
      \psfrag{lba'}{$l_{ba'}$}
      \psfrag{lbb'}{$l_{bb'}$}
      \psfrag{lbc'}{$l_{bc'}$}
      \psfrag{ap}{$s_{a'}$}  
      \psfrag{bp}{$s_{b'}$}
      \psfrag{cp}{$s_{c'}$}
      \psfrag{gap}{$g_{a'}$} 
      \psfrag{gbp}{$g_{b'}$}
      \psfrag{gcp}{$g_{c'}$}
      \psfrag{lpaap}{$l'_{aa'}$}
      \psfrag{lpabp}{$l'_{ab'}$}
      \psfrag{lpacp}{$l'_{ac'}$}
      \psfrag{lpcap}{$l'_{ca'}$}
      \psfrag{lpcbp}{$l'_{bb'}$}
      \psfrag{lpccp}{$l'_{bc'}$}
      \psfrag{lpbap}{$l'_{ba'}$}
      \psfrag{lpbbp}{$l'_{cb'}$}
      \psfrag{lpbcp}{$l'_{cc'}$}
      \psfrag{E(D)}{$E(\Delta)$}
      \psfrag{E(Dp)}{$E(\Delta')$}
      \begin{minipage}{.55\textwidth}
        \includegraphics[width=7cm]{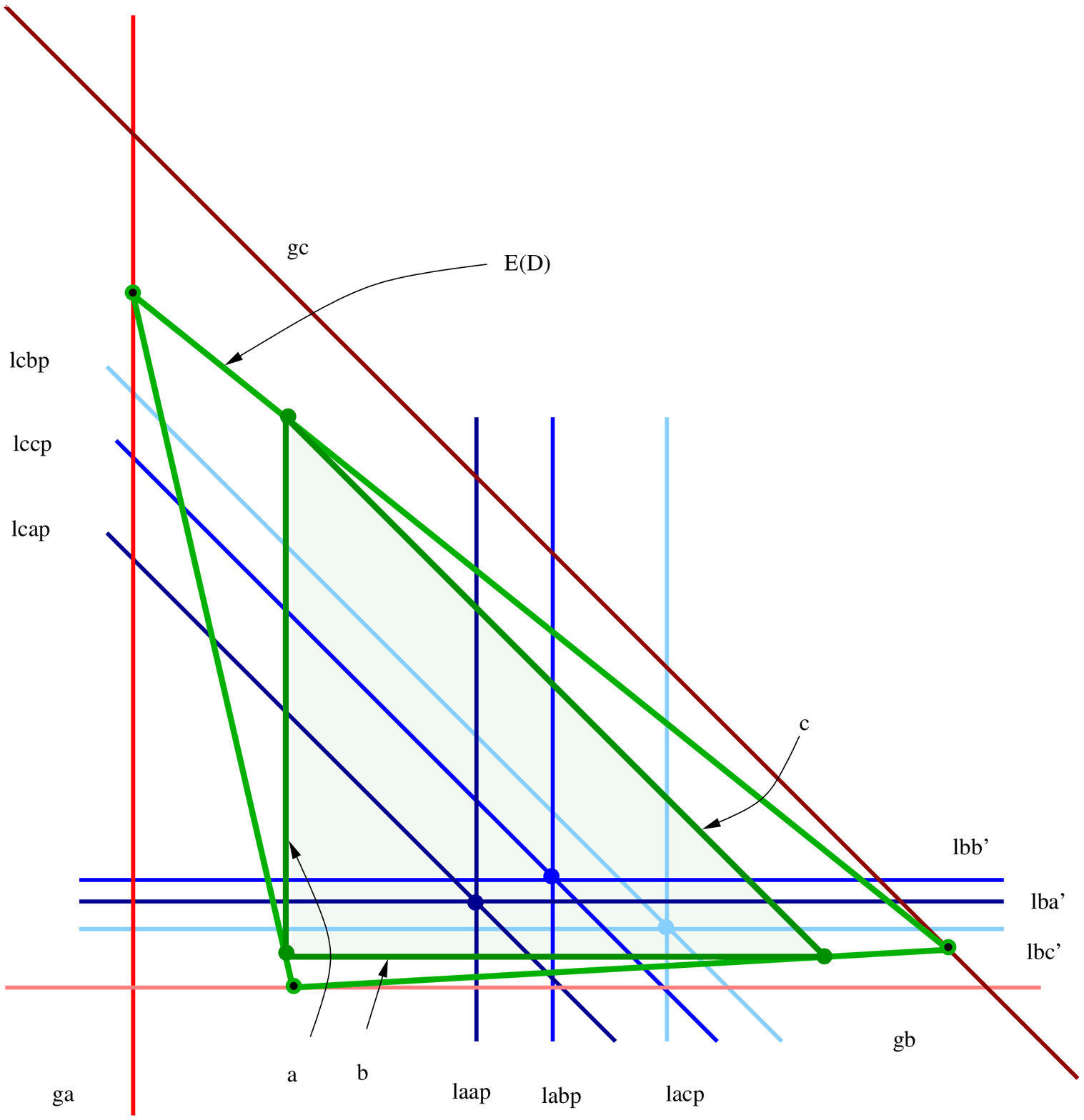}
      \end{minipage}
      \hspace{-.12\textwidth}
      \begin{minipage}{.55\textwidth}

        \vspace*{-.05\textheight}
        \includegraphics[width=7cm]{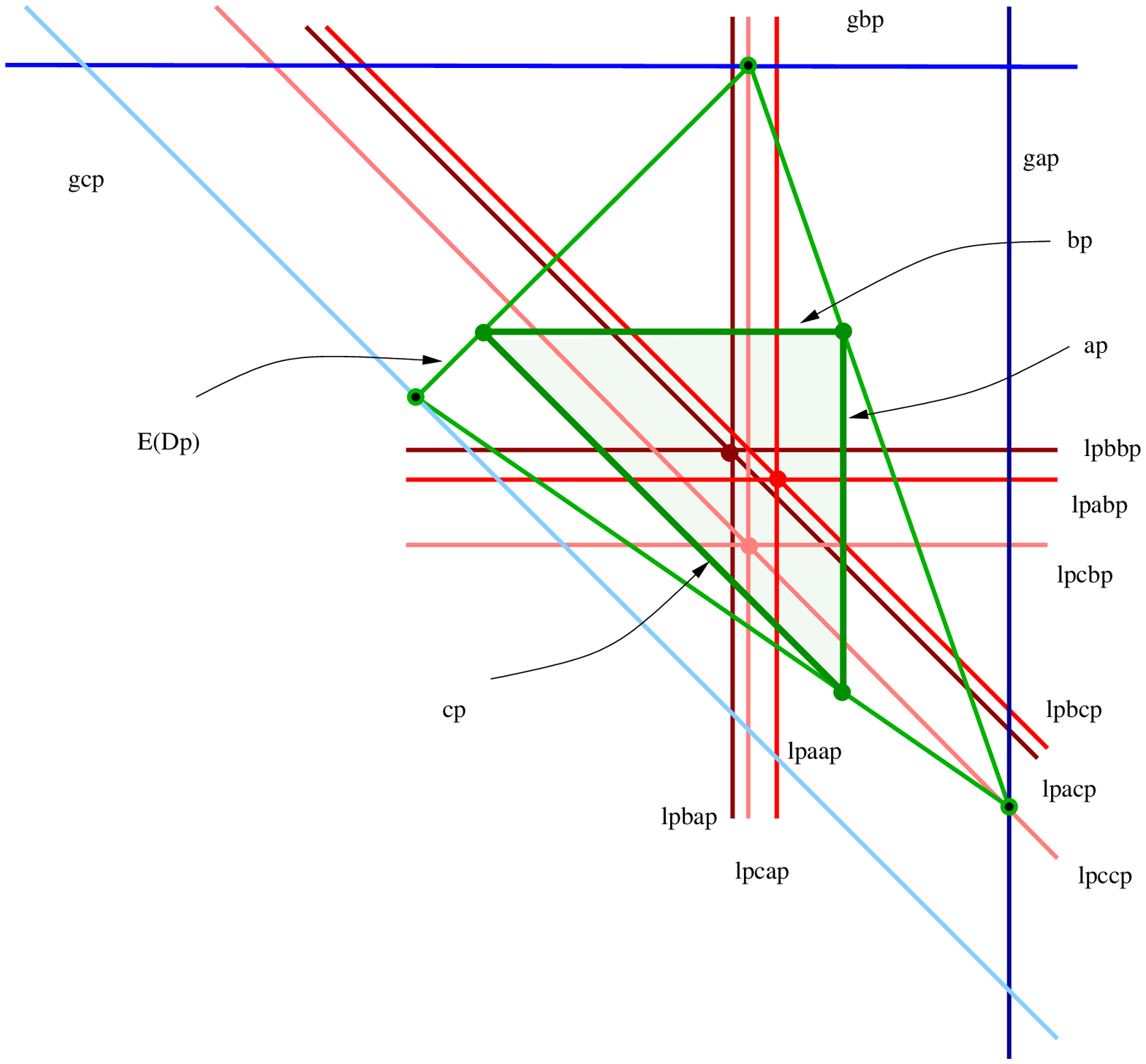}
      \end{minipage}
      \caption{Construction of the triangles}
      \label{fig:triangle-ratios}
    \end{minipage}
  \end{figure}
  In the following let  the index $x$  always run  through $\{a,b,c\}$
  and and $y$ through $\{a',b',c'\}$.  Fix  two triangles $\Delta$ and
  $\Delta'$ and   let $s_a, s_b, s_c$  be  the sides   of $\Delta$ and
  $s_{a'}, s_{b'}, s_{c'}$ the sides  of $\Delta'$.  By a  translation
  in each of the two factors we can assume that  they both contain the
  origin.  Denote the nine  ratios by $r_{xy}$ for $x\in\{a,b,c\}$ and
  $y\in\{a',  b', c'\}$.   See  Figure \ref{fig:triangle-ratios}.   To
  simplify    the   notation     we     introduce    the    parameters
  $R_{xy}:=\frac{r_{xy}}{1-r_{xy}}$.
  
  Let $g_x$ be a line outside $\Delta$ parallel to $s_x$ at a distance
  $\delta_x$.  These three lines will  afterwards contain the vertices
  of  $E(\Delta)$, which   is a triangle  containing   the vertices of
  $\Delta$ in  its   edges.  Similarly, define    the line $g_{y}$  at
  distance $\delta_y$ from $s_y$ for $\Delta'$.
  
  Let $l_{ay}$ define a line parallel to $s_a$ lying on the other side
  of $a$ as $g_a$ at distance $R_{ay}\delta_a$  from $s_a$.  Similarly
  define the lines   $l_{by}$ and $l_{cy}$  parallel  to $b$ and  $c$.
  Thus,   any segment  starting on  $g_x$  and ending  on  $l_{xy}$ is
  divided   by  $s_x$ with  a   ratio of $r_{xy}$.   For  the triangle
  $\Delta'$ we define lines $l'_{ay}$  at distance $1/R_{ay}$ parallel
  to $s_y$ and on  the other side  as $g_y$.  Finally, define (outward
  pointing) normal vectors $n_x$ and $n_y$ and levels $\lambda_x$, and
  $\lambda_{y}$    such     that    points   $u\in    s_x$     satisfy
  $\skp{n_x}{u}-\lambda_x=0$   and    points   $v\in   s_y$    satisfy
  $\skp{n_y}{v}-\lambda_y=0$.
  
  Consider now e.g.   the ratio $r_{ab'}$.  Choose  a vertex  $v_a$ of
  $E(\Delta)$ on $g_a$,  a point $w_a$  on the line  $l'_{ab'}$ and in
  the interior of $\Delta'$, a  vertex $v'_{b'}$ of $E(\Delta')$ lying
  on $g_{b'}$ and a point $w'_{b'}$ in the interior of $\Delta$ on the
  line $l_{ab'}$.  The  points  $w_a$  and $w'_{b'}$ will  become  the
  corresponding points to $v_a$ and $v'_{b'}$ under the maps $\beta_0$
  and $\beta_1$  of (\ref{conditionB}).     The part of    the segment
  $h_{ab'}$   between $v_a$ and $w_{b'}'$    lying inside $\Delta$ has
  length   $r_{ab'}|h_{ab'}|$  (where $|h|$ denotes   the  length of a
  segment $h$) and the part of the segment $h_{b'a}$ between $v'_{b'}$
  and $w_a$  inside  $\Delta'$ has length  $(1-r_{ab'})|h_{b'a}|$.  So
  the condition set by  the ratio $r_{ab'}$  will be satisfied by this
  choice of $w_a$ and $w'_{b'}$.
  
  To satisfy all conditions on the ratios  that involve $w_a$, we have
  to choose $w_a$ such that it lies as well on  the lines $l_{ay}$ and
  in  the interior of $\Delta$.   Similar  conditions hold for the two
  other  points  inside $\Delta$  and   for the  three  points  inside
  $\Delta'$.   Therefore,   finding a   feasible  solution  amounts to
  finding a solution to the following set of $18$ linear equations and
  six linear inequalities.
  \begin{align*}
    \lambda_x&=\skp{n_x}{w'_{y}}+R_{xy}\delta_x\\
    \lambda_{y}&=\skp{n_{y}}{w_{x}}+1/R_{xy}\delta_{y}\\
    0&<\delta_x, \delta_y
  \end{align*}
  for  all    $x\in\{a,b,c\}$    and $y\in\{a',b',c'\}$.    Here   the
  coordinates of     the  points  $w_x$, $w'_y$  and     the distances
  $\delta_x$,  $\delta_y$ are the free  variables,  and the ratios are
  the parameters.  The   first and the   second set  of  equations are
  connected via the ratios.   As the equations and inequalities depend
  smoothly  on the nine parameters,  it suffices for  the proof of the
  theorem to show that there exists at  least one feasible solution of
  this   system.     Such  a     solution   is  shown    in     Figure
  \ref{fig:triangle-ratios}  and  in  Table~\ref{tab:feasiblesolution}
  (for  some  fixed  product  of    two triangles,   but this can   be
  projectively transformed to any other).
  
  To finally  obtain $E(\Delta)$ we   have to  choose vertices on  the
  lines  $g_a$,  $g_b$,  and $g_c$  such  that the  edges  contain the
  vertices   of   $\Delta$.    Unless  the  distances    $\delta_{a}$,
  $\delta_{b}$, and $\delta_{c}$ are too large compared to the size of
  $\Delta$, there are always two solutions to this problem (one of the
  solutions   for  the  $E_{33}$     mentioned   above is  given    in
  Table~\ref{tab:feasiblesolution},   the     other   is  obtained  by
  reflection),   which    depend    continuously   on   the  distances
  $\delta_x,\delta_y$ (There is no  solution otherwise).  Similarly we
  can construct $E(\Delta')$.
\end{proof} 

\begin{table}[b]
  \scriptsize
  \begin{minipage}[b]{.4\textwidth}
    \centering
    \begin{alignat*}{5} 
      [&& 1 &&\;\;  1 &&\;\;  1 &&\;\;  1]\\[-.06cm]
      [&& 1 &&\;\;  1 &&\;\;  1 &&\;\; -1]\\[-.06cm]
      [&& 1 &&\;\;  1 &&\;\; -1 &&\;\;  1]\\[-.06cm]
      [&& 1 &&\;\;  1 &&\;\; -1 &&\;\; -1]\\[-.06cm]
      [&& 1 &&\;\; -1 &&\;\;  1 &&\;\;  1]\\[-.06cm]
      [&& 1 &&\;\; -1 &&\;\;  1 &&\;\; -1]\\[-.06cm]
      [&& 1 &&\;\; -1 &&\;\; -1  &&\;\; 1]\\[-.06cm]
      [&& 1 &&\;\; -1 &&\;\; -1 &&\;\; -1]\\[-.06cm]
      [&&-1 &&\;\;  1 &&\;\;  1 &&\;\;  1]\\[-.06cm]
      [&&-1 &&\;\;  1 &&\;\;  1 &&\;\; -1]\\[-.06cm]
      [&&-1 &&\;\;  1 &&\;\; -1 &&\;\;  1]\\[-.06cm]
      [&&-1 &&\;\;  1 &&\;\; -1 &&\;\; -1]\\[-.06cm]
      [&&-1 &&\;\; -1 &&\;\;  1 &&\;\;  1]\\[-.06cm]
      [&&-1 &&\;\; -1 &&\;\;  1 &&\;\; -1]\\[-.06cm]
      [&&-1 &&\;\; -1 &&\;\; -1 &&\;\;  1]\\[-.06cm]
      [&&-1 &&\;\; -1 &&\;\; -1 &&\;\; -1]\\[-.06cm]
      [&& 3/5 &&\;\; 9/5 &&\;\;-3/5 &&\;\;-3/5]\\[-.06cm]
      [&& 9/5 &&\;\;-3/5 &&\;\;-3/5
      &&\;\; 3/5]\\[-.06cm]
      [&&-3/5 &&\;\;-9/5 &&\;\; 3/5 &&\;\;3/5]\\[-.06cm]
      [&&-9/5 &&\;\; 3/5 &&\;\; 3/5 &&\;\;-3/5]\\[-.06cm]
      [&&-3/5 &&\;\; 3/5 &&\;\; 3/5 &&\;\;9/5]\\[-.06cm]
      [&&-3/5 &&\;\;-3/5 &&\;\;-9/5 &&\;\;3/5]\\[-.06cm]
      [&& 3/5 &&\;\;-3/5 &&\;\;-3/5 &&\;\;-9/5]\\[-.06cm]
      [&& 3/5 &&\;\; 3/5 &&\;\; 9/5 &&\;\;-3/5]
    \end{alignat*}
    \caption{An $E_{44}$   from regular squares,  but  not  satisfying
      (\ref{ass:ratioequal}) of Theorem~\ref{thm:symmetry}.}
    \label{tab:less_symmetry}
  \end{minipage}
  \hspace{.03\textwidth}
  \begin{minipage}[b]{.58\textwidth}
    \centering
    \begin{alignat*}{5}
      [&& 1&&\;\; 0&&\;\; 1&&\;\; 0]\\
      [&& 1&&\;\; 0&&\;\; 0&&\;\; 0]\\
      [&& 1&&\;\; 0&&\;\; 0&&\;\; 1]\\
      [&& 0&&\;\; 0&&\;\; 1&&\;\; 0]\\
      [&& 0&&\;\; 0&&\;\; 0&&\;\; 0]\\
      [&& 0&&\;\; 0&&\;\; 0&&\;\; 1]\\
      [&& 0&&\;\; 1&&\;\; 1&&\;\; 0]\\
      [&& 0&&\;\; 1&&\;\; 0&&\;\; 0]\\
      [&& 0&&\;\; 1&&\;\; 0&&\;\; 1]\\
      [&&    9/247&&\;\; 289/247&&\;\;  819/1387&&\;\;  364/1387]\\
      [&&  289/247&&\;\;- 51/247&&\;\;  364/1387&&\;\;  204/1387]\\
      [&& - 51/247&&\;\;   9/247&&\;\;  204/1387&&\;\;  819/1387]\\
      [&&  153/494&&\;\;  34/247&&\;\;  169/1387&&\;\; 1764/1387]\\
      [&&   34/247&&\;\;  21/38&&\;\;  1764/1387&&\;\; -546/1387]\\
      [&&   21/38&&\;\; 153/494&&\;\;  -546/1387&&\;\;  169/1387]
    \end{alignat*}
    \caption{The coordinates  of  a feasible  non-degenerate solution.
      See Figure~\ref{fig:ProdTriangles}  for a   drawing of  the  two
      factors.}
    \label{tab:feasiblesolution}
  \end{minipage}
\end{table}
From this construction  method the proof of Theorem~\ref{thm:relspace}
is straightforward:

\begin{proof}[Proof of Theorem~\ref{thm:relspace}]
  All triangles  in $\R^2$ are  projectively equivalent. Therefore, up
  to projective equivalence, there is  only one geometric  realisation
  of an orthogonal product of two triangles. Thus, in the following we
  can fix our preferred orthogonal  product of two triangles and count
  the degrees of freedom  for  adding the remaining vertices   without
  having to worry about projective equivalence anymore.  But according
  to Theorem~\ref{thm:triangleratios} we  have, for any  choice of two
  triangles, nine degrees of freedom  for the choice of the  remaining
  vertices.
\end{proof}

\begin{remark}
  There     might  still be    geometric   realisations  of a polytope
  combinatorially equivalent   to $E_{33}$ that are   not projectively
  equivalent to a polytope   containing an orthogonal product   of two
  triangles.  Thus, a priori Theorem~\ref{thm:relspace} describes only
  a subset of the whole realisation space $\rs_{proj}(E_{33})$.
\end{remark}

\subsubsection{The 24-cell.}

From  our method  to   realise the   $E$-construction of  products  of
polygons we can obtain new geometric realisations of the $24$-cell.

\paragraph{A 4-parameter family of 24-cells.}

For $m,n>3$  we cannot determine  the degrees of  freedom in the above
way anymore. Taking the $mn$ ratios as input we obtain $2mn$ equations
and   $m+n$  inequalities for  only   $3(m+n)$ variables.  This is not
merely  a   problem of  the method.  There    are  in  fact additional
restrictions on a realisation, as the  lengths of the segments from an
interior  point  to the   vertices  of the $E$-construction  cannot be
viewed as independent variables anymore (consider e.g.\ a square and a
pair of opposite vertices of its $E$-construction).  However, also for
the $24$-cell it    is   not  difficult to  construct     projectively
non-equivalent geometric realisations.

\begin{figure}[t] 
\includegraphics[height=.24\textheight]{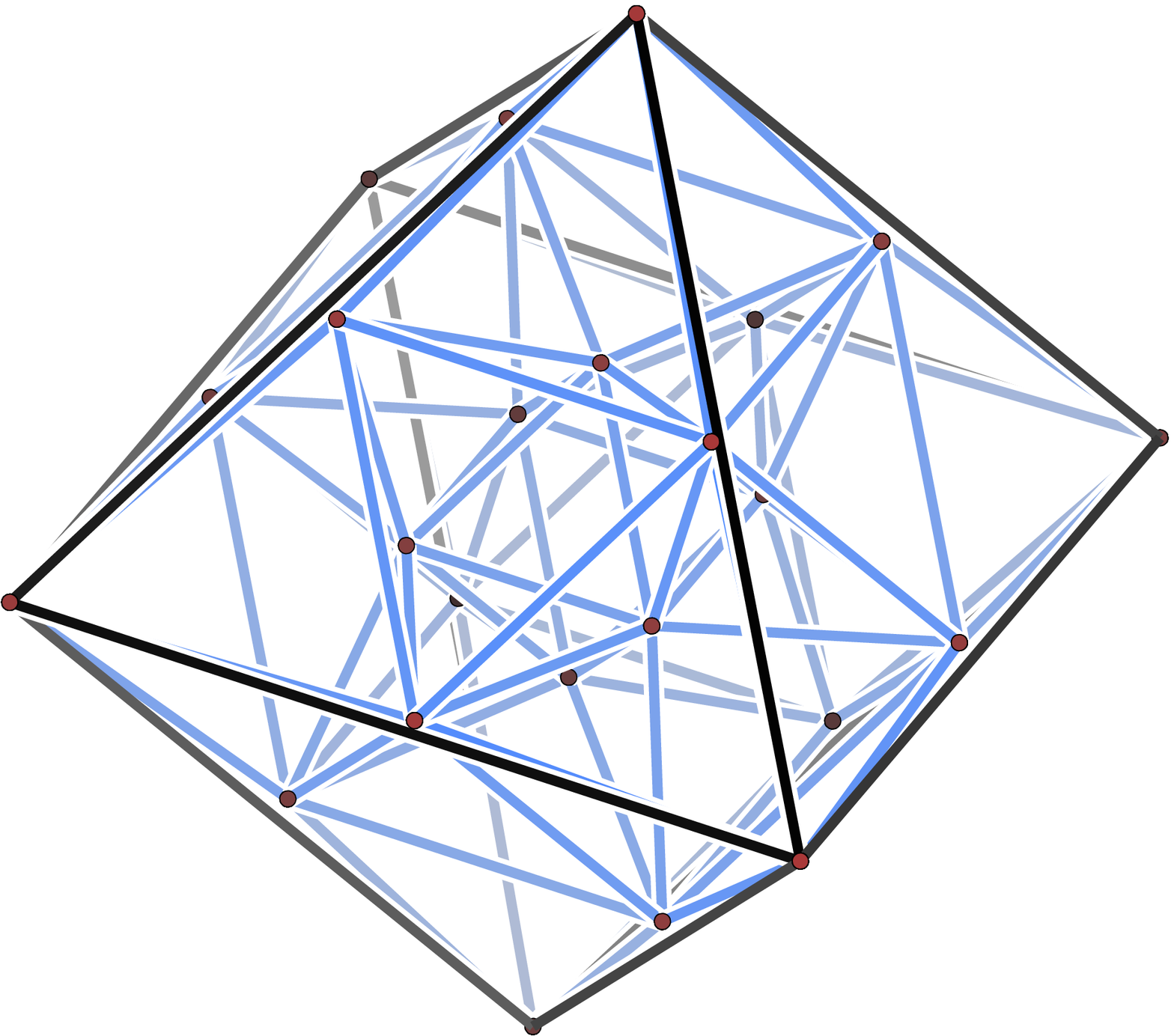}
\hspace{-.03\textwidth}
\includegraphics[height=.24\textheight]{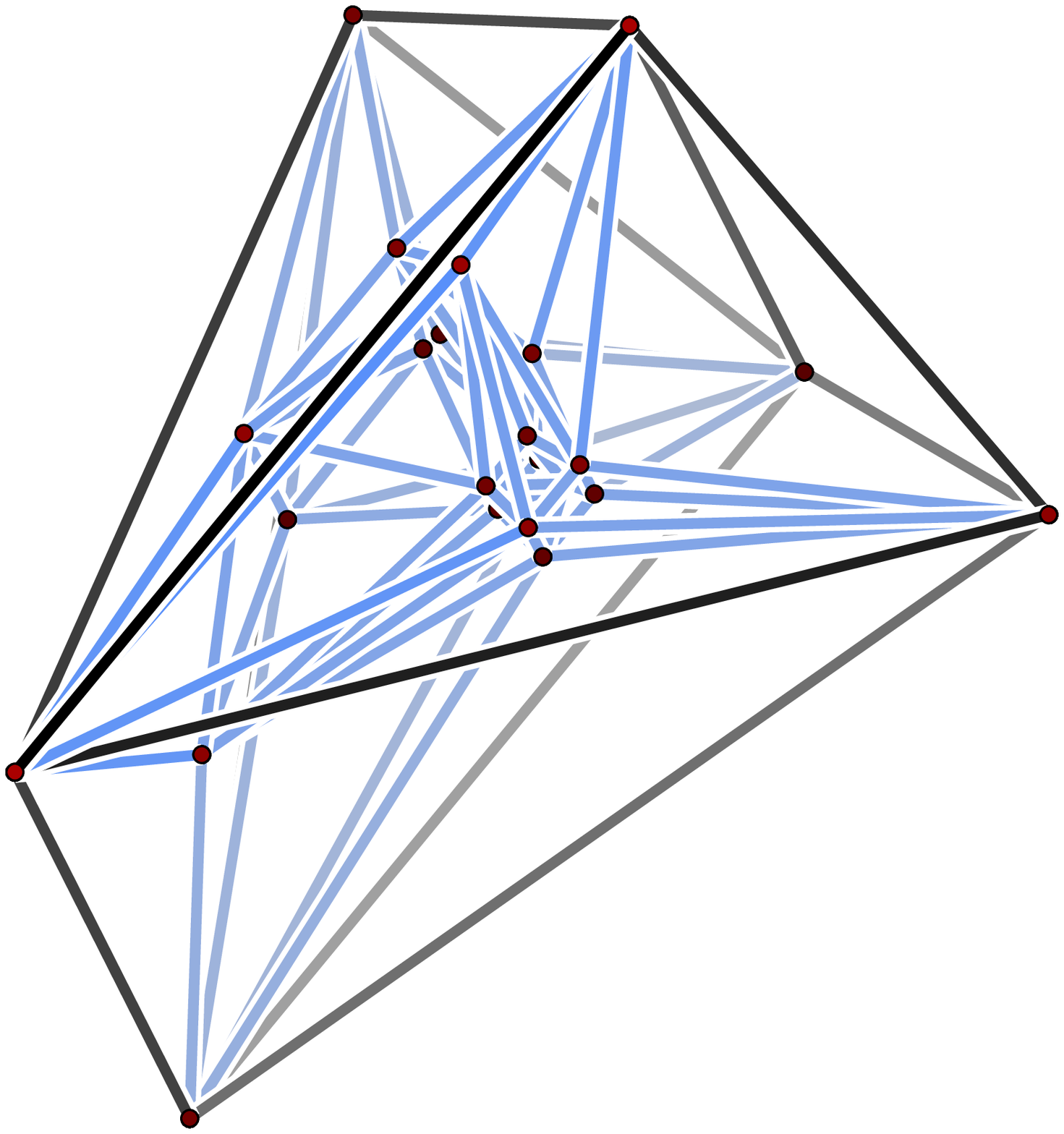}
\hspace{-.5cm}
\caption{The regular 24-cell and  one other realisation
  from the family of Table~\ref{fig:Family24Cells}.}
\label{fig:24Schlegel}\label{fig:nop}
\end{figure}

Table~\ref{fig:Family24Cells}    shows  a   simple   example    of   a
$4$-parameter  family of  $24$-cells, where all  four parameters range
between $-1$ and  $1$.  This family spans  a $4$-dimensional subset of
the projective realisation  space, which can  be seen in the following
way.  The vertex set of the  regular $24$-cell contains the vertex set
of  three different standard cubes: If  you set  all parameters to $0$
then (in the order  given in Table~\ref{fig:Family24Cells}) the  first
sixteen, the last sixteen and the first  and last eighth vertices each
form  a standard  cube.   Their  $2$-faces (squares)  are not  anymore
present  in     the   $24$-cell,  but     they   still    lie on     a
codimension-2-subspace,   which  is    preserved   by  any  projective
transformation   (e.g.      vertices     $15,  16,      17, 18$     in
Table~\ref{fig:Family24Cells}).  Letting  the parameters diverge  from
$0$ destroys some of these ``internal'' squares, necessarily resulting
in projectively different $24$-cells.  This  can also  be seen in  the
Schlegel diagrams  in Figure~\ref{fig:24Schlegel}:  Observe the  three
squares  contained in  the octahedral  face  on which the  polytope is
projected.
\begin{table}[b]
  \scriptsize
  \begin{minipage}[t]{.55\textwidth}
    \begin{alignat*}{5}
      [&&        -1&&        -1&&        -1&&        -1]\\[-.05cm]
      [&&         1&&         1&&        -1&&        -1]\\[-.05cm]
      [&&         1&&        -1&&         1&&        -1]\\[-.05cm]
      [&&        -1&&         1&&         1&&        -1]\\[-.05cm]
      [&&         1&&        -1&&        -1&&         1]\\[-.05cm]
      [&&        -1&&         1&&        -1&&         1]\\[-.05cm]
      [&&        -1&&        -1&&         1&&         1]\\[-.05cm]
      [&&         1&&         1&&         1&&         1]\\[-.05cm]
      [&&         1&&        -1&&        -1&&        -1]\\[-.05cm]
      [&&        -1&&         1&&        -1&&        -1]\\[-.05cm]
      [&&        -1&&        -1&&         1&&        -1]\\[-.05cm]
      [&&        -1&&        -1&&        -1&&         1]\\[-.05cm]
      [&&         1&&         1&&         1&&        -1]\\[-.05cm]
      [&&         1&&         1&&        -1&&         1]\\[-.05cm]
      [&&         1&&        -1&&         1&&         1]\\[-.05cm]
      [&&        -1&&         1&&         1&&         1]\\[-.05cm]
      [&&       a_1&&       b_1&&       a_2&&\;\;-2-b_2]\\[-.05cm]
      [&&       a_1&&       b_1&&     2-a_2&&       b_2]\\[-.05cm]
      [&&       a_1&&       b_1&&       a_2&&     2-b_2]\\[-.05cm]
      [&&       a_1&&       b_1&&\;\;-2-a_2&&       b_2]\\[-.05cm]
      [&&       a_1&&     2-b_1&&       a_2&&       b_2]\\[-.05cm]
      [&&\;\;-2-a_1&&       b_1&&       a_2&&       b_2]\\[-.05cm]
      [&&       a_1&&\;\;-2-b_1&&       a_2&&       b_2]\\[-.05cm]
      [&&     2-a_1&&       b_1&&       a_2&&       b_2]
    \end{alignat*}
    \caption{Vertices of a family of $24$-cells}
    \label{fig:Family24Cells}
  \end{minipage}
  \hspace{.05\textwidth}
  \begin{minipage}[t]{.4\textwidth}
    \begin{alignat*}{5}
      [&&  -1&&   5/4&&  -1&&   1],\\[-.05cm]
      [&&  -1&&   5/4&&  -1&&  -1],\\[-.05cm]
      [&&  -1&&   5/4&&   1&&  -1],\\[-.05cm]
      [&&  -1&&   5/4&& 5/3&&   1],\\[-.05cm]
      [&&  -1&&    -1&&  -1&&   1],\\[-.05cm]
      [&&  -1&&    -1&&  -1&&  -1],\\[-.05cm]
      [&&  -1&&    -1&&   1&&  -1],\\[-.05cm]
      [&&  -1&&    -1&& 5/3&&   1],\\[-.05cm]
      [&&   1&&    -1&&  -1&&   1],\\[-.05cm]
      [&&   1&&    -1&&  -1&&  -1],\\[-.05cm]
      [&&   1&&    -1&&   1&&  -1],\\[-.05cm]
      [&&   1&&    -1&& 5/3&&   1],\\[-.05cm]
      [&&   1&&\;\; 23/12&& -1&& 1],\\[-.05cm]
      [&&   1&& 23/12&&  -1&&  -1],\\[-.05cm]
      [&&   1&& 23/12&&   1&&  -1],\\[-.05cm]
      [&&   1&& 23/12&& 5/3&&   1],\\[-.05cm]
      [&&-1/2&&  -1/2&&-3/2&& 1/2],\\[-.05cm]
      [&&-1/2&&  -1/2&&\;\;-5/6&&\;\; -3/2],\\[-.05cm]
      [&&-1/2&&  -1/2&&17/6&&-1/2],\\[-.05cm]
      [&&-1/2&&  -1/2&& 1/2&& 5/2],\\[-.05cm]
      [&&-3/2&&  -5/6&&-1/2&&-1/2],\\[-.05cm]
      [&& 1/2&&  -3/2&&-1/2&&-1/2],\\[-.05cm]
      [&& 5/2&&   1/2&&-1/2&&-1/2],\\[-.05cm]
      [&&-1/2&&  10/3&&-1/2&&-1/2],
\end{alignat*}
\caption{{A $24$-cell without any projective automorphisms.}}
\label{fig:NoProjAutos}
\end{minipage}
\end{table} 

\begin{remark}
  Clearly, not all  possible     realisations of the  $24$-cell    are
  contained  in  this    $4$-parameter  family.    The $24$-cell    in
  Table~\ref{fig:NoProjAutos} is also a result of the construction and
  has no projective automorphisms.
\end{remark}

\subsection{Fatness of polytopes}

The classification of $f$- and flag vectors for polytopes in dimension
$d\ge 4$ is   an important unsolved problem in   polytope theory.  See
\cite{Bayer87},  \cite{Ziegler02} for some  background on this problem
and overviews of the known results.

Ziegler \cite{Ziegler02}  proposed to look  at the  following quantity
(called the ``fatness'' of a polytope $P$) on the entries of these two
vectors.
\begin{align*}
F(P):=\frac{f_1+f_2-20}{f_0-f_3-10}
\end{align*}
where  $(f_0,f_1,f_2,f_3)$ is   the  $f$-vector  of  any  $4$-polytope
different  from  the  simplex (in  \cite{EKZ02}  there  is  a slightly
different definition).  The fatness  of  polytopes produced  from  the
$E$-construction applied  to simple  $4$-polytopes is  bounded by~$6$,
cf.  \cite[p.~3]{EKZ02}.  Eppstein,  Kuperberg, and Ziegler provided a
polytope $Q$ resulting  from a gluing of  $600$-cells that has fatness
around   $5.073$  in  the     definition  of  \cite{Ziegler02}    (The
$E$-construction also works  for some  non-simple polytopes,  but  all
known examples don't have a  higher  fatness).  They also showed  that
for regular CW  $3$-spheres  fatness is unbounded.   But  they neither
found polytopes with fatness higher than $5.073$ nor an upper bound on
fatness for arbitrary polytopes.

For our family $E_{mn}$ we get according to the $f$-vector computation
in (\ref{eq:fvector}):
\begin{align*}
  F(E(C_m\times C_n))=\frac{12mn-20}{2mn+2m+2n-10}\longrightarrow 6
\end{align*}
for $m,n\rightarrow\infty$.  Thus  for  $m,n\ge 10$ our polytopes  are
``fatter'' than  the   above mentioned example  from  \cite{EKZ02}. As
products of polygons  are simple,  our family  of  polytopes is ``best
possible'' within   this setting.   However,  Ziegler~\cite{Ziegler04}
recently has constructed a  class of polytopes (not ``$E$-polytopes'')
with  fatness arbitrarily close to   $9$ by considering projections of
products of polygons to $\R^4$.

\renewcommand{\refname}{References}
 
{\footnotesize
\providecommand{\bysame}{\leavevmode\hbox to3em{\hrulefill}\thinspace}

\bibliographystyle{amsalpha}

}
\end{document}